\documentclass[12pt]{article}
\usepackage{amsfonts,amssymb,latexsym}
\normalbaselines
\begin{document}
\newcommand{\me}{\mod(q-\eps)}
\newcommand{\NN}{{\cal N}}
\newcommand{\ZZ}{{\cal Z}}
\newcommand{\FF}{{\cal F}}
\newcommand{\OO}{{\cal O}}
\newcommand{\LL}{{\cal L}}
\newcommand{\PP}{{\cal P}}
\newcommand{\DD}{{\cal D}}
\newcommand{\Aa}{{\cal A}}
\newcommand{\Bb}{{\cal B}}
\newcommand{\II}{{\cal I}}
\newcommand{\JJ}{{\cal J}}
\newcommand{\QQ}{{\cal Q}}
\newcommand{\MM}{{\cal M}}
\newcommand{\Cb}{\Bbb C}
\newcommand{\Sb}{{\Bbb S}}
\newcommand{\Qb}{{\Bbb Q}}
\newcommand{\Ob}{\Bbb O}
\newcommand{\eps}{\varepsilon}
\newcommand{\Zb}{\Bbb Z}
\newcommand{\Nb}{\Bbb N}
\newcommand{\Mb}{\Bbb M}
\newcommand{\de}{\delta}
\newcommand{\Goldie}{\mbox{Goldie}}
\newcommand{\Center}{\mbox{Center}}
\newcommand{\Ad}{\mbox{Ad}}
\newcommand{\ad}{\mbox{ad}}
\newcommand{\Sp}{\mbox{Spec}}
\newcommand{\Fract}{\mbox{Fract}}
\newcommand{\Rad}{\mbox{Rad}}
\newcommand{\mod}{\mbox{mod}}
\newcommand{\ux}{{\underline{x}}}
\newcommand{\uz}{{\underline{z}}}
\newcommand{\uy}{{\underline{y}}}
\newcommand{\ub}{{\underline{b}}}
\newcommand{\ua}{{\underline{a}}}
\newcommand{\um}{{\underline{m}}}
\newcommand{\un}{{\underline{n}}}
\newcommand{\ta}{\tilde{a}}
\newcommand{\tb}{\tilde{b}}
\newcommand{\tc}{\tilde{c}}
\newcommand{\td}{\tilde{d}}

\newcommand{\oDD}{\overline{\DD}}
\newcommand{\Ms}{{\mbox{Maxspec}}}
\newcommand{\rk}{{\mbox{rank}}}
\newcommand{\Disc}{\mbox{Disc}}
\newcommand{\Ipd}{I(\pi)_\DD}
\newcommand{\Ihpd}{\widehat{I}(\pi)_\DD}
\newcommand{\Rc}{R_\chi}
\newcommand{\Rcd}{R_{\chi,\DD}}
\newcommand{\Rpd}{R_{\pi,\DD}}
\newcommand{\Bpd}{B_{\pi,\DD}}
\newcommand{\Rme}{R_{\mu\eps}}
\newcommand{\Sme}{S_{\mu\eps}}
\newcommand{\Epd}{E_{\pi\DD}}
\newcommand{\Bme}{B_{\mu\eps}}
\newcommand{\Zme}{Z_{\mu\eps}}
\newcommand{\Pme}{\PP_{\mu\eps}}
\author{A.~N.~Panov.}
\date{\it{Mathematical Department, Samara State University,\\
 ul.Akad.Pavlova 1, Samara, 443011, Russia\\
panov@info.ssu.samara.ru}}
\title{
Quantum solvable algebras. Ideals and representations at roots of 1}
\maketitle
\footnote{The work is supported by RFFI grant  02-01-00017.}

{\bf Abstract}. Let $R$ be a quantum solvable algebra.
It is proved that every prime ideal $I$ which is stable with respect to
quantum adjoint action is completely prime, and $\Fract(R/I)$ 
is isomorphic to the skew field of fractions of an 
algebra of twisted polynomials.
We study the correspondence between symplectic leaves  and irreducible 
representations. A Conjecture of De Concini-Kac-Procesi on the dimension 
of irreducible representations is proved for 
sufficiently great $l$. 

\section {Introduction.}
A quantum solvable algebra is an iterated $q$-skew extension,
constructed by derivations and diagonal automorphisms. 
Among examples are Quantum Matrices,
Quantum Weyl algebra, Quantum Heisenberg algebra, the positive part
$U_q({\bf n})$ of the quantum group $U_q({\bf g})$ 
and also some their subalgebras
and factor algebras. 
Representations of these algebras were treated separately
[DCKP1],[DC-P],[JZ].
The main approach is a correspondence between irreducible representations
of a quantum algebra 
and symplectic leaves of Poisson manifold of the center.
 
In the paper we study quantum solvable algebras in general setting
(Definition 3.1), 
obeying some natural and easily checkable conditions (Conditions 3.2-3.4).
We consider the specialisations of these algebras at primitive
$l^{\rm{th}}$ of unity $\eps$. 

We make some assumptions on $l$.
First, we require that 
the elements $x_1^l,\ldots,x_{n+m}^l$ lie in the center of $R_\eps$
and $\eps$ is rather great:  $\eps$ is a 
point of good reduction
of stratification of prime spectra for intederminate $q$ 
(Definition 3.10). 
Second, we require that $l$ is 
relatively prime with all minors of integer matrix
$\Sb$ (Definition 3.1) and with all integers $s_1,\ldots,s_n$ (Condition 3.2). 
Almost any $\eps$ that obeys the second conditions is a 
point of good reduction (Theorem 3.15). 
We say that $\eps$ is admissible, if it obeys the above conditions
(Definition 3.20). 

We study prime ideals stable with respect to
quantum adjoint action (Definition 2.1). 
We study their intersections with
components of filtration of quantum solvable algebras
(Theorem 2.12).      
All prime ideals which are stable with respect to quantum adjoint action
are completely prime (Theorem 3.24).
We get a stratification of these ideals.
This helps us to study irreducible representations 
and symplectic leaves (Theorem 5.3).
In particular, we prove that the dimension of any irreducible
representation equals to $l^{d/2}$ where $d$ is the dimension of 
the corresponding
symplectic leaf. It was conjectured by C.De Concini, V.G.Kac, and 
C.Procesi
in the papers [DCKP1,4.5] and [DC-P,25.1].

We are very grateful for the referees for numerous pertinent remarks
and advices. 
 
\section{On quantum adjoint action and stable ideals}
In this section we study skew polynomial extensions.
Let $R$ be a ring, $\tau$ be an endomorphism of $R$, and $\delta$ be
a $\tau$-derivation of $R$ (i.e. $\delta(ab)=\delta(a)b+\tau(a)\delta(b)$
for all $a,b\in R$).
A skew polynomial extension of a ring $R$ is the ring  
$E=R[\ux;\tau,\delta]$ generated by indeterminate $\ux$ and $R$ obeying
$\ux a=\tau(a)\ux+\delta(a)$ for all $a\in R$ (see [MC-R,2.1], [GL1]).
The main statement of this section is Theorem 2.12 on the reduction
of prime ideals in skew extensions.

Throughout the paper we use the following notations.
Let $K$ be a field of zero characteristic and let $q$ be an indeterminate.
We denote by $C$ a localization of $K[q,q^{-1}]$ over some finitely
generated denominator set and $F=\Fract(C)$. 
Let $\Gamma$ be the cyclic subgroup 
$\{q^n\}_{n\in\Zb}$ in the group of invertible elements $C^*$ of $C$.
Let $\eps$ be a primitive $l^{\rm{th}}$ root of unity.
Denote by $K(\eps)$ the field extension of $K$ 
(if $\eps\in K$, we put $K(\eps)=K$). 
Suppose that $C$ admits specialisation via
$q\mapsto\eps$ in the field $K(\eps)$.

Throughout this section $R$ is an unitary $C$-algebra and 
a free $C$-module.
We consider specialisation
$\rho_\eps:R\to R_\eps=R \bmod(q-\eps)$.

Further we differentiate between 
the elements in $R$ with their images in $R_\eps$, 
using the following notations.
For $x\in R$ we denote $x_\eps=x\bmod(q-\eps)R$ or
for $y\in R_\eps$ we denote by $\underline{y}$ an element in 
its preimage in $R$.

Let $u$ be an element in $R$ such that $u_\eps=u\bmod(q-\eps)$
lies in the center $Z_\eps$ of $R_\eps$.
For any $\ua\in R$ the element $\ad_u(\ua)=u\ua-\ua u$ lies in $(q-\eps)R$.
Whence $\frac{\ad_u(\ua)}{(q-\eps)}\in R$.
We denote 
$$
\DD_u(a)=\frac{\ad_u(\ua)}{(q-\eps)}\bmod(q-\eps)\eqno(2.1)
$$
for $a=\ua\bmod(q-\eps)\in R_\eps$. It is easy to check that the definition
(2.1) is independent of the choice of preimage $\ua$ of $a$.\\
{\bf Definition 2.1}. We say that $\DD_u:R_\eps\to R_\eps$ is the quantum
adjoint action of $u\in R$.\\
{\bf Remark}. The derivations $\DD_u$ and their properties
were studed in the papers
[DCKP1-2] and [DC-P]. There is no 
definition of quantum adjoint action in the above
papers. The term of quantum coadjoint action were used for the action of
 the subgroup ${\mbox exp}(t\DD_u)$ on maximal spectrum of $Z_\eps$.

The following properties are straightforward to check.\\
{\bf Property 2.2}. $\DD_u$ is a derivation of $R_\eps$. That is
$\DD_u(ab)=\DD_u(a)b+a\DD_u(b)$ for all $a,b\in R_\eps$.\\
{\bf Property 2.3}. If $u\bmod(q-\eps)=u_1\bmod(q-\eps)$, then
$u_1-u=(q-\eps)r$, $r\in R$ and $\DD_{u_1}=\DD_u+\ad_r$.
This implies that an ideal $\PP$ is $\DD_{u_1}$-stable whenever it is
$\DD_u$-stable.\\
{\bf Property 2.4}. If $u_\eps,v_\eps\in Z_\eps$ for $u,v\in R$, then
$$
\DD_{uv}(a)=\DD_u(a)v_\eps+u_\eps\DD_v(a)
$$
{\bf Property 2.5}. The bracket 
$$
\{u_\eps,v_\eps\}=\DD_u(v_\eps)=-\DD_v(u_\eps)
$$
is a Poisson bracket on $Z_\eps$.

Consider a skew polynomial extension
$E=R[\ux;\tau,\delta]$ where $\tau$ is an automorphism and $\delta$ is
$\tau$-derivation of $R$. 
We assume that $\tau\vert_C=\mbox{id}$, $\delta\vert_C=0$
and suppose that $E$ is a $q^s$-skew extension of $R$. That is
$\delta\tau=q^s\tau\delta$ for some $s\in \Zb$ ([G],[GL1-2]).
One can extend the automorphism
$\tau$ to an automorphism of $E$, putting $\tau(\ux)=q^{-s}\ux$. 
We save the notation $\tau$ for this extension.
Notice that $E_\eps=E\me=R_\eps[x;\tau,\delta]$ where $x=\ux\me$\\
{\bf Assumptions 2.6.} 
1) As above $\eps$ is a primitive $l^{\rm{th}}$ root of unity.
We suppose that $GCD(s,l)=1$ if $s\ne 0$.\\
2) Suppose that $x^l$, $x=\ux\bmod(q-\eps)$ lie in the center 
 of $E_\eps$. We shall denote by $\DD_x$ the derivation
$\DD_{\ux^l}$ defined by (2.1).\\
3) We assume that $\tau$ is a $\Gamma$-diagonalizable automorphism
of $R$.
That is, there exists a $C$-basis $\{a_\alpha\}$ of $R$ such that 
$\tau(a_\alpha) = q^{m_\alpha} a_\alpha$ for some $m_\alpha\in \Zb$.
It follows that the similar basis also exists in $E$.

 One can extend $\tau$ and $\de$ to an automorphism 
and $\tau$-derivation of $R_\eps$. We save the same notations for these 
extensions.
Further we use usual notations \\
1) $(n)=(n)_{q^s}=\frac{q^{sn}-1}{q^s-1}$ (if $s=0$, put $(n)=n$);\\
2) $(n)!=(n)_{q^s}!=(1)(2)\cdots(n)$;\\
3) ${n\choose k}={n\choose k}_{q^s}=\frac{(n)!}{(k)!(n-k)!}$.\\
{\bf Lemma 2.7}. Let $\eps,R$ be as in 
Assumptions 2.6(1,2); then $\delta^l=0\bmod(q-\eps)$,
$\tau^l=\mbox{id}\ \bmod(q-\eps)$.\\
{\bf Proof.} 
One can prove ([G,6.2],[GL,2.5]) for all $a\in R$
$$
\ux^la =\sum_{i=0}^l{l\choose i}\tau^{l-i}\delta^i(a)\ux^{l-i}=
\tau^l(a)\ux^l+\sum_{i=1}^{l-1}{l\choose i}\tau^{l-i}\delta^i(a)\ux^{l-i}+
\delta^l(a).\eqno(2.2)
$$
Since the element $x^l=\ux^l\bmod(q-\eps)$ lies in the center
$Z_\eps$, then 
$\tau^l(a)=a\bmod(q-\eps)$ and $\delta^l(a)=0\bmod(q-\eps)$.$\Box$

We shall use the following notations 
$$\Delta = \frac{\de^l}{q-\eps}:R_\eps\to R_\eps,$$
$$\theta = \frac{\tau^l-\mbox{id}}{q-\eps} :E_\eps\to E_\eps.$$ 
Remark that $\Delta$  is a linear operator on $R_\eps$ and 
$\theta$ is a derivation of $E_\eps$.
For any integer $m$ we denote
$$\um=\left(\frac{q^{ml}-1}{q-\eps}\right)_{q=\eps}=ml\eps^{-1}.$$
Notice that, if $\tau(a)=q^ma$ for some $a\in E$ and $m\in\Zb$, then
$\theta(a_\eps)=\um a_\eps$.

Consider the localization $R'$ of $R$ over the set 
$$\{f(q)\in C: f(\eps)\ne 0\}.$$
{\bf Lemma 2.8}. Let $\eps, R$ be as in Assumptions 2.6(1,2);
then $\frac{\delta^k}{(k)!}(R')\subset R'$ for any positive integer $k$.\\
{\bf Proof.} For $s=0$ the statement is trivial. Let $s\ne 0$.
 Let $k$ be a positive integer.
Suppose that $k=\xi l+\eta$ for integers $\xi$ and $\eta$ with $\xi$ positive
and $0\le\eta<l$. By above Lemma 2.7, $\delta^l(R)\subset (q-\eps)R$.
For $a$ in $R$ we have 
$$\de^k(a)=\de^{\xi l+\eta}(a)=\de^\eta(\de^l)^\xi(a)\subset
(q-\eps)^\xi\de^\eta(R)\subset(q-\eps)^\xi R.$$
On the other hand, 
$(q-\eps)$ divides $(i)=\frac{q^{is}-1}{q^s-1}$ if and only if
$l$ divides $i$. It follows $(k)!=(q-\eps)^\xi c(q)$ where $c(\eps)\ne 0$.
This implies the assertion.$\Box$

For all $g\in R$ we denote  
$$\Pi_n(g)=\tau^{n-1}(g)\cdots\tau(g) g.$$
{\bf Lemma 2.9.} Let $\eps, R$ obey Assumptions 2.6(1,2).
Let $a\in R$ and let $I_a$ be the minimal $\tau$-ideal generated by 
$a$ in $R'$. 
Let $n,m$ be positive integers.
We assert that  \\
$$\frac{\delta^n(a^m)}{(n)!}=\left\{\begin{array}{cc}
\Pi_n(\delta(a))\bmod I_a & \mbox{for}\quad m=n,\\
0\ \bmod I_a&\mbox{for}\quad n<m.\end{array}\right.
$$
{\bf Proof.} We shall prove the lemma by induction on $n$.
The case $n=1$ is easy.
Suppose that statement is true for all positive integers less then $n$.
We are going to prove for $n$.
We need the formula [G,6.2]:
$$
\delta^n(ab)=\sum_{i=0}^n{n\choose i}\tau^{n-i}\delta^i(a)\delta^{n-i}(b).
\eqno(2.3)
$$
We get
$$\delta^n(a^m)=\delta(a^{m-1}a)=
\sum_{i=0}^n{n\choose i}\tau^{n-i}\delta^i(a^{m-1})\delta^{n-i}(a).$$
We obtain
$$
\frac{\delta^n(a^m)}{(n)!}=
\sum_{i=0}^{n-2}\tau^{n-i}\frac{\delta^i(a^{m-1})}{(i)!}
\frac{\delta^{n-i}(a)}{(n-i)!}
+ \frac{\tau\delta^{n-1}(a^{m-1})}{(n-1)!}\delta(a)+
\frac{\delta^n(a^{m-1})}{(n)!}a.\eqno (2.4)
$$
Let $n<m$. By the induction hypothesis,
$\frac{\delta^i(a^{m-1})}{(i)!}\in I_a$ for $i\le n-1$. 
By Lemma 2.8, $\frac{\delta^{n-i}(a)}{(n-i)!}\in R'$.
This proves that all terms in (2.4), apart from the last one, lie in $I_a$.
As for the last term, it obviously lies in $I_a$. This concludes the 
case $n<m$.

Let $m=n$.
Similarly to the previous case, we get
$$
\frac{\delta^n(a^{n})}{(n)!}= 
 \frac{\tau\delta^{n-1}(a^{n-1})}{(n-1)!}\delta(a)\bmod I_a =
\tau(\Pi_{n-1}(\delta(a)))\delta(a)\bmod I_a=
\Pi_n(\delta(a))\bmod I_a.
$$
$\Box$

Let $I$ be an ideal in $E_\eps$. Denote $J=I\bigcap R_\eps$.\\
{\bf Lemma 2.10}. 
Let Assumptions 2.6(1-3) hold.
Suppose that $I$ is $\theta$-stable.\\
1) $I$ is $\tau$-stable.\\ 
2) The ideal $J$ is $(\tau,\delta)$-stable. 
If in addition $I$ is $\DD_x$-stable, then
$J$ is a $\Delta$-stable ideal in $R_\eps$.\\
{\bf Proof.}\\
1) The automorphism $\tau: E\to E$ is diagonalizable.
Let $\{a_\alpha\}$ be a basis such that each $a_\alpha$ is $\tau$-eigenvector:
$\tau(a_\alpha)=q^{m_\alpha}a_\alpha$.
Consider the reduction of $C$-basis ${a_\alpha}$ modulo $q-\eps$.
We get a $\Cb$-basis $a_{\alpha\eps}=a_\alpha\bmod(q-\eps)$.
The derivation $\theta$ is diagonalizable
and $\theta(a_{\alpha\eps})=\um_\alpha a_{\alpha\eps}$.
Every $\theta$-eigenvector lies in a span of basis elements $a_{\alpha\eps}$
corresponding 
to a common $\theta$-eigevalue $\um$. Then all these $a_{\alpha\eps}$
have the common $\tau$-eigenvalue $\eps^m$. This proves that
every $\theta$-eigenvector is also $\tau$-eigenvector
(the converse is not true).    

Let $a$ be an element in $I$. 
Decompose $a$
into a sum of $\theta$-eigenvectors with different eigenvalues
(which we call $\theta$-components of $a$).
Since $I$ is $\theta$-stable, all $\theta$-components also lie in $I$.
We saw above that this $\theta$-components are $\tau$-eigenvectors.
These prove that $\tau(a)\in I$. We conclude that $I$ is $\tau$-stable.\\
2) Since $I$ is $(\tau,\theta)$-stable, the ideal $J$ is also 
$(\tau,\theta)$-stable. Let $a$ be in $J$.
Since $xa=\tau(a)x+\de(a)$, then 
$xa, \tau(a)x\in I$. Whence $\delta(a)\in J$.
This proves that $J$ is $\delta$-stable.

Suppose that $I$ is $\DD_x$-stable. 
By direct calculations,
$$\DD_x(a)=\theta(a)x^l+
\sum_{i=1}^ {l-1}\frac{{l\choose i}_{q^s}}{q-\eps}\tau^{l-i}(\delta^i(a))
x^{l-i} + \Delta(a).\eqno(2.5)
$$
Since $\theta(a)\in J$ and $\DD_x(a)\in I$, then $\Delta(a)\in J$.
$\Box$\\
{\bf Lemma 2.11.} Let $\tau,\delta$ be as above; then 
$\delta(\theta(a))=(\theta+\underline{s})\delta(a)$ for all 
$a\in R_\eps$.\\
{\bf Proof}. Putting $q=\eps$ in 
$$\de\frac{\tau^l-\mbox{id}}{q-\eps}=
 \frac{q^{sl}\tau^l-\mbox{id}}{q-\eps}\de=
\left(q^{sl}\frac{\tau^l-\mbox{id}}{q-\eps}+\frac{q^{sl}-1}{q-\eps}\right)\de,
$$
we obtain $\de\theta=(\theta+\underline{s})\de$.
$\Box$\\ 
{\bf Theorem 2.12}. Let $R$ be a Noetherian $C$-algebra and a free 
$C$-module. Suppose that $E=R[\ux;\tau,\delta]$
is a $q^s$-skew extension.\\
1) Let $R$ obey Assumption 2.6(3).
Let $\II$ be a prime $\tau$-stable ideal in $E$. 
Then the ideal
$\JJ:=\II\bigcap R$ is also prime.\\
2)  
Let $\eps$ and $ R$ obey Assumptions 2.6(1-3).
Let $I$ be a prime $(\theta,\DD_x)$-stable  
ideal in $E_\eps$. Then
the ideal $J= I\bigcap R_\eps$ is also prime.\\
{\bf Proof.}\\
1) Since $\II$ is $\tau$-stable, then $\JJ$ is $(\tau,\delta)$-stable. 
The ideal $\II$ is prime $\tau$-stable;
 then $\JJ$ is $(\tau,\delta)$-prime.
By [G,6.5], $\JJ$ is $\tau$-prime. It follows that 
$\JJ$ is semiprime and minimal prime ideals over $\JJ$ form a $\tau$-orbit.
Since $\JJ$ is $\tau$-stable, $\JJ$ is prime. \\  
2) The proof divides in two points 2i) and 2ii). \\ 
2i) In this point we are going to prove that $J$ is semiprime.
By Lemma 2.10, 
the ideal $J$ is a $(\tau,\delta)$-stable.
Then $EJ$ is a two-sided ideal in $E_\eps$.
 The ideal $J$ is a $(\tau,\delta)$-prime ideal in $R_\eps$
[GL2, 2.1(vi)]. 

Let $N$ be an ideal in $R_\eps$ such that $N/J=\mbox{Rad}(R_\eps/J)$.
Hence, $N^n\subseteq J$ for some positive integer $n$.
We are going to prove that $N=J$. 

The radical $N$ is stable under actions of automorphisms and derivations
[MC-R,14.2.3]. In particular
$N$ is $(\tau,\theta)$-stable.

Suppose that $N$ is $\delta$-stable. 
Then $E_\eps N$ is a two-sided ideal in $E_\eps$.
We see that $(E_\eps N)^n\subseteq E_\eps N^n\subseteq E_\eps J\subseteq I$.
Since $I$ is prime, we have $E_\eps N\subseteq I$. 
Notice that $N=(E_\eps N)\bigcap R_\eps$ and $J=I\bigcap R_\eps$.
Whence $N\subseteq J$.
This proves that $J$ is semiprime.

The above observation shows that 
it remains to prove that $N$ is $\delta$-stable.
Consider the image $\LL=\delta(N)\bmod(N)$ in $R_\eps/N$
and suppose that $\LL$ is non zero.
Since
$$
\delta(N)=\delta(R_\eps NR_\eps)=
\delta(R_\eps)N R_\eps+\tau(R_\eps)\delta(N)R_\eps+
\tau(R_\eps)\tau(N) R_\eps,
$$ 
we see that $\LL$ is an ideal in $B=R_\eps/N$. 
Moreover, by Lemma 2.11,
$\theta\delta(N)=\delta(\theta(N))-\underline{s}\delta(N)$.
Whence
the ideal $\LL$ is $(\tau,\theta)$-stable.
The algebra $B$ is semiprime. Its
Goldie quotient algebra $A=\Goldie(B)$ is a direct sum of matrix algebras
over division rings $A=\oplus_{i=1}^NA_i$. The 
automorphism $\tau$ acts on
components $A_i$ by permutations.
The ideal $A\LL A$ is the sum $A'=\oplus_{i=1}^tA_i$ which is $\tau$-stable.
 There exist the elements 
$r_1,\ldots,r_k, p_1,\ldots,p_k\in A$ and $h_1,\ldots,h_k\in \LL$ such that
$$
1_{A'}=r_1h_1p_1+\cdots r_kh_kp_k.
$$
Multiplying this equality by suitable
regular elements of $B$, we deduce 
that $\LL$ contains a regular element of $A'$,
 say $h$. This implies that $\Pi_n(h)\ne 0$ with $n$ as above.

Denote by $a$ an element of $R_\eps$ such that $a\in N$, $\delta(a)=h$.
Denote by $\underline{a}\in R$ a representative of $a\in R_\eps$.
By Lemma 2.9, we have 
$$\frac{\delta^n(\underline{a}^n)}{(n)!}=
\Pi_n(h)+\underline{b}$$
with some $\underline{b}\in I_{\underline{a}}.$ 
As in Lemma 2.8 decompose $n=\xi l+\eta$, $0\le \eta<l$ and $\xi$,$\eta$ are 
positive integers. 
Denote 
$$\underline{d}=\frac{\delta^n(\underline{a}^n)}{(q-\eps)^\xi}=c(q)(
\Pi_n(h)+\underline{b})$$
with $c(\eps)\ne 0$.
Take the above equality modulo $q-\eps$.
Notice that $R'\bmod(q-\eps)=R\bmod(q-\eps)=R_\eps$.
The element $d=\underline{d}\bmod(q-\eps)$
equals to $A'$-regular element $c(\eps)\Pi_n(h)$ modulo $N$.

On the other hand,
$$ \underline{d}=\frac{\delta^n(\underline{a}^n)}{(q-\eps)^\xi}=
\left(\frac{\delta^l}{q-\eps)}\right)^\xi \delta^\eta(\underline{a}^n)=
\Delta^\xi\delta^\eta(a^n)
\bmod(q-\eps),$$
$$d=\Delta^\xi\delta^\eta(a^n).$$
Recall that $a^n\in N^n\subseteq J$ and $J$ is $(\delta,\Delta)$-stable,
by Lemma 2.10(2).
It follows that $d\in J$. Whence $d=0\bmod N$.
A contradiction.
The ideal $N$ is $\delta$-stable and $J$ is semiprime.\\
2ii) In this point we shall prove that $J$ is prime.

Since $J$ is semiprime, one can present it as 
an intersection of prime ideals 
$J=Q_1\bigcap\cdots\bigcap Q_m$. 
If a derivation preserves $J$, then it preserves all $Q_i$ (see 
[MC-R,14.2.3] and [D,3.3.2]).
This implies that all ideals $Q_i$ are $\theta$. 
Whence $Q_i$ is $\tau$-stable, as in the proof of Lemma 2.10(1).
Retain the notations of the point 2i) (with now $N=J$).
 Decompose
$A=\oplus_{i=1}^mA_i$, $A_i=Ae_i$.
 Automorphism $\tau$-acts on the system of central minimal idempotents
$e_1,\ldots,e_m$ by permutations.
Since $\tau(A_i)\subset A_i$, then $\tau(e_i)=e_i$ and
$$\delta(e_i)=\delta(e_i^2)=\delta(e_i)e_i+\tau(e_i)\delta(e_i)=
2e_i\delta(e_i).
$$  
Since $K$ has characteristic zero, this proves $\delta(e_i)=0$.
Whence $\delta(Q_i)\subset Q_i$ for all $i$.
Since $J$ is $(\tau,\delta)$-prime, $J=Q_i$ . Finally, 
$J$ is prime.\\
$\Box$

\section{Stratification of  $\DD$-stable ideals}

In this section we are going to get stratification of prime  $\DD$-stable 
ideals at roots of unity. 
This provides a reduction to the case of 
algebra of twisted Laurent polynomials.
We prove that all prime  $\DD$-stable ideals are completely prime.
 
Recall $K$ is a field of zero characteristic, $q$ is an indeterminate
and $C$ is a localization of $K[q,q^{-1}]$ over some finitely 
generated denominator set.
 Let $\Sb=(s_{ij})$ be a $M\times M$ integer skew-symmetric matrix.
Denote $q_{ij}=q^{s_{ij}}$ and form the matrix $\Qb=(q_{ij})$.\\
{\bf Definition 3.1}([P1],[P2]).
We say that a ring $R$ is quantum solvable over $C$, if
$R$ is generated by the elements 
$x_1, x_2,\ldots, x_n, x_{n+1}^{\pm 1},\ldots,
x_{n+m}^{\pm 1}$ with $M=n+m$
such that the monomials $x_1^{t_1}\cdots x_n^{t_n}x_{n+1}^{t_{n+1}}\cdots
x_{n+m}^{t_{n+m}}$ with $t_1,\ldots,t_n\in\Nb$ and 
$t_{n+1},\ldots,t_{n+m}\in\Zb$ 
form a free $C$-basis with the 
relations\\
1) $x_ix_j=q_{ij}x_jx_i$ for all $i$ and $n+1\le j\le n+m$ ;\\
2) for $1\le i<j\le n$ the formula holds 
$$x_ix_j=q_{ij}x_jx_i+r_{ij},\eqno (3.1)$$
where $r_{ij}$ is an element of the subalgebra
$R_{i+1}$ generated by $x_{i+1},\ldots, x_n, x_{n+1}^{\pm 1},\ldots,
x_{n+m}^{\pm 1}.$

There exists a chain of subalgebras
$R=R_1\supset R_2\supset\cdots\supset R_n\supset 
R_{n+1}$.
The last algebra $R_{n+1}$ is an algebra of twisted Laurent polynomials
generated by $x_{n+1}^{\pm 1},\ldots,x_m^{\pm 1}$.
One can prove that a quantum 
solvable algebra is an iterated skew extension of $C$
(see [MC-R],[GL1]).
This means that  
 for all $i\in [1,n]$ the map 
$\tau_i: x_j\mapsto q_{ij}x_j, i>j$ is extended to an automorphism of
$R_{i+1}$ and the map $\delta_i:x_j\mapsto r_{ij}$ is extended to 
a $\tau_i$-derivation of $R_{i+1}$. 
Each algebra $R_i$ is the skew extension $R_{i+1}[x_i;\tau_i.\delta_i]$.
All automorphisms $\tau_i$ are identical on $C$ and all 
$\tau_i$-derivations $\delta_i$ are equal to zero on $C$.
A quantum solvable algebra is a Noetherian domain [MC-R,1.2.9].
We put some more conditions on $R$. This conditions are comparable with 
more general conditions Q1-Q4 of [P2].\\
{\bf Condition 3.2}. We require that $R$ is an iterated $q$-skew extension
in terms of [GL1-2].
This means that $\delta_i\tau_i=q_i\tau_i\delta_i$ for some
$q_i=q^{s_i}$, $s_i\in\Zb$.\\  
{\bf Condition 3.3}. All automorphisms $\tau_i$ are extended to diagonal
automorphisms of $R$ (i.e the monomials 
$x_1^{t_1}\cdots x_n^{t_n}x_{n+1}^{s_1}\cdots
x_{n+m}^{s_m}$, $t_i\in\Nb$, $s_i\in\Zb$ are eigenvectors) with 
eigenvalues in $\Gamma=\{q^k: n\in\Zb\}$. 
Denote by $H$ the group of diagonal automorphisms
of $R$ generated by $\tau_1,\ldots,\tau_{n+m}$. We shall refer
a common eigenvector for all $\tau_1,\ldots,\tau_{n+m}$ as a 
$H$-weight vector.
\\
{\bf Condition 3.4}.  
Let $\Lambda$ be the subset
of the set of all roots of unity in the algebraic closure 
$\overline{K}$ such that
$\eps\in\Lambda$ if $\eps$ is a primitive $l^{\rm{th}}$ root of unity
and the elements  
$x_1^l,\ldots,x_n^l$ modulo $q-\eps$ lie in the center $Z_\eps$ of $R_\eps$.
We require that $\Lambda$ is infinite.

The last Condition claims that $R_\eps$ is finite over its center for 
an  infinite set of specialisations of $C$. 
Such algebras $R$ are called pure $C$-algebras
in the paper [P1]. 
The Conditions 3.2-3.4 imply the Conditions Q1-Q4 of the paper [P2]
(in the case $C$ is a localization of $K[q,q^{-1}]$ over some finitely
generated denominator set). See [P2] for examples. One gets  
the following stratification of prime ideals 
with zero intersection with $C$. \\
{\bf Definition 3.5}. We say that two elements $a,b$ $q$-commute if 
$ab=q^kba$ for some integer $k$.\\
{\bf Theorem 3.6}[P2, Theorem 3.4]. Let $R$  be a quantum solvable algebra
obeying Conditions 3.2-3.4.\\
1) There exists a finite set $\Mb=\{\PP_\mu\}$
of semiprime $H$-stable ideals with $\PP\bigcap C=0$ for all $\mu$, 
and denominator subsets
$S_\mu\subset R/\PP_\mu$ with $S_\mu$ generated by $k+m$ $q$-commuting 
$H$-weight elements.
Here,
 $\mu=(i_1,\ldots,i_{k+1})$ where
 $k$ and $\{i_s: 1\le s\le k+1\}$ are non-negative integers, and
  $k+i_1+\cdots+i_{k+1}=n$.
The localization 
$(R/\PP_\mu)S_\mu^{-1}$ is isomorphic to a factor algebra of an algebra
of twisted Laurent polynomials.\\
2) Any prime ideal $\II$ of $R$ with zero intersection with $C$ 
(recall that indeed all such ideals are completely prime)
contains a unique ideal $\PP_\mu$ such that 
$S_\mu\bigcap(\II/\PP_\mu)=\emptyset$. 

This yields a stratification of $\Sp(R)$.
We recall here some of the details, referring to 
[P2, Theorem 3.4] for complete proof.

Fix any integer $i_1$ with $0\le i_1\le n$. 
Consider the subset  $L_1=\{x_n,x_{n-1},\ldots,x_{n-i_1+1}\}$
of $i_1$ generators of $R$.
Let $\JJ_1$ be the minimal semiprime ideal containing $L_1$.
Denote by 
$X_1$ the set of minimal prime ideals over $\JJ_1$.
We decompose $X_1=X_1'\bigcup X_1''\bigcup X_1'''$
where
$$X_1'=\{\QQ\in X_1:\QQ\bigcap C= 0, x_{n-i_1}\notin\QQ, \},$$ 
$$X_1''=\{\QQ\in X_1:\QQ \bigcap C\ne 0, x_{n-i_1}\notin\QQ\},$$
$$X_1'''=\{\QQ\in X_1: x_{n-i_1}\in\QQ\}.$$ 
We denote
$$\PP_{i_1}=\bigcap_{Q\in X_1'}Q.\eqno(3.2)$$ 
If $X_1'=\emptyset$, then there is no prime ideal $\II$ of $R$ such
that  $\II\supset L_1$, $\II\bigcap C=0$, and $x_{n-i_1}\notin \II$. 
We add $x_{n-i_1}$ in $L_1$
and begin the stratification process from the beginning.
 
The ideal $\PP_{i_1}$ is $H$-stable, contains $L_1$, does not contain 
$x_{n-i_1}$ and $\PP_{i_1}\bigcap C=0$.
Denote by $S_{1*}$ the denominator subset $\{x_{n-i_1}^m\}_{m\in\Nb}$
in $R/\PP_{i_1}$.
Consider the localizaton $R^{(1)}=(R/\PP_{i_1})S_{1*}^{-1}$.
The algebra $R^{(1)}$ is generated by 
$x_1,\ldots,x_{n-i_1-1},x_{n-i_1}^{\pm 1},R_{n+1}$ and admits the filtration
$$
R^{(1)}=R_1^{(1)}\supseteq\cdots\supseteq
R_{n-i_1}^{(1)}\supseteq R_{n-i_1+1}^{(1)}=\cdots=R_n^{(1)}=R_{n+1}^{(1)}
$$
with
$R_j^{(1)}= (R_j/{\PP_{i_1}\bigcap R_j})S_{1*}^{-1}$ for $j\le n-i_1$
and $R_{n+1}^{(1)}=R_{n+1}/\PP_{i_1}\bigcap R_{n+1}$.
The subalgebra $R_j^{(1)}$ is generated by $x_j$ and $R_{j+1}^{(1)}$.

There exists a localization of $C$ over some finitely generated denominator 
subset $N_1\subset C$ such that the algebra $R^{(1)}N_1^{-1}$ over 
$C_1=CN_1^{-1}$ has the
new system of $H$-weight generators
$x_1',\cdots,x'_{n-i_1-1},x_{n-i_1}^{\pm 1},
x_{n+1}^{\pm 1},\ldots , x_{n+m}^{\pm 1}$ with the relations
$x'_jx_{n-i_1}=q_{j,n-i_1}x_{n-i_1}x'_j$ [P2, Lemma 3.3].
Here $q_{j,n-i_1}$ is the $(j,n-i_1)$-th entry of the matrix $\Qb$
(see Definition 3.1) 
The algebra $R_j^{(1)}N_1^{-1}$ is generated by $x_j'$ and 
$R_{j+1}^{(1)}N_1^{-1}$.
Denote by $S_1$ the denominator set  generated by $N_1$ and $S_{1*}$, 
and retain the notation $R^{(1)}$ for the localization $R^{(1)}N_1^{-1}$.
We get $R^{(1)}=(R/\PP_{i_1})S_1^{-1}$.

We continue the stratification process.
Now let $i_2$ be any integer with $0\le i_2\le n-i_1-1$.
Consider the subset
$L_2=\{x'_{n-i_1-1},\ldots,x'_{n-i_1-i_2}\}$ of $i_2$ generators of 
$R^{(1)}$.
Let $\JJ_2$ be the minimal semiprime ideal of $R^{(1)}$ containing $L_2$.
We denote by $X_2$ the set of minimal prime ideals
 over $\JJ_2$.
As above we decompose
$X_2=X_2'\bigcup X_2''\bigcup X_2'''$
where
$$X_2'=\{\QQ'\in X_2:\QQ'\bigcap C_1= 0, x'_{n-i_1-i_2-1}\notin\QQ, \},$$ 
$$X_2''=\{\QQ'\in X_2:\QQ' \bigcap C_1\ne 0, x'_{n-i_1-i_2-1}\notin\QQ\},$$
$$X_2'''=\{\QQ'\in X_2: x'_{n-i_1-i_2-1}\in\QQ\}.$$
We denote
$$\PP'_{i_1i_2}=\bigcap_{Q'\in X_2'}Q'.\eqno (3.3)$$
There exists a semiprime $H$-stable ideal
$\PP_{i_1i_2} $ is $R$ such that
$$(\PP_{i_1i_2}/{\PP_{i_1}})S_1^{-1}=\PP_{i_1i_2}'.$$
Preserve the notation $S_1$ for the image of $S_1$ in $R/\PP_{i_1i_2}$.
Consider the denominator set 
$S'_{2*}$ generated by $x'_{n-i_1-i_2-1}$ in $(R/\PP_{i_1i_2})S_1^{-1}$.
The elements from $S_1$ $q$-commute with $S'_{2*}$.
There exists a positive integer $t$ such that 
$x'_{n-i_1-i_2-1}x_{n-i_1}^t\in R/\PP_{i_1i_2}$.
Denote by  
$S_{2*}$ the denominator subset generated by $S_1$ and 
$x'_{n-i_1-i_2-1}x_{n-i_1}^t$
in $R/\PP_{i_1i_2}$. 
Denote
$R^{(2)}=(R/\PP_{i_1i_2})S_1^{-1}(S_{2*}')^{-1}=
(R/\PP_{i_1i_2})S_{2*}^{-1}$.
The algebra 
$R^{(2)}$ is generated by
$$
x_1',\cdots,x'_{n-i_1-i_2-2},(x'_{n-i_1-i_2-1})^{\pm 1}, x_{n-i_1}^{\pm 1},
x_{n+1}^{\pm 1},\ldots , x_{n+m}^{\pm 1}$$ 
and admits the filtration
$$
R^{(2)}=R_1^{(2)}\supseteq\cdots\supseteq
R_{n-i_1-i_2}^{(2)}=\cdots=R_{n-i_1}^{(2)}\supseteq
R_{n-i_1+1}^{(2)}=\cdots=R_n^{(2)}=R_{n+1}^{(2)}.
$$
There exists 
localization of $C$ over some finitely generated denominator 
subset $N_2\subset C$ such that the algebra $R^{(2)}N_1^{-1}$  has
the new system of $H$-weight generators
$$x_1'',\cdots,x''_{n-i_1-i_2-2},(x'_{n-i_1-i_2-1})^{\pm 1}, 
x_{n-i_1}^{\pm 1},x_{n+1}^{\pm 1},\ldots , x_{n+m}^{\pm 1}$$
where
$x''_j$ $q$-commutes with $x'_{n-i_1-i_2-1}$ and $x_{n-i_1}$.
As above we denote
by $S_2$ the denominator set generated by $S_{2*}$ and $N_2$.
We retain the notation 
 $R^{(2)}$ for localization $R^{(2)}N_2^{-1}$.
 We get
 $R^{(2)}=(R/\PP_{i_1i_2})S_{2}^{-1}$.

We continue the stratification process. Finally, we get a
semiprime $H$-stable ideal
$\PP_\mu=\PP_{i_1i_2\ldots i_{k+1}}$ and a denominator subset
$S_\mu$ in $R/\PP_\mu$. 
The subset $S_\mu$ is generated by $q$-commuting $H$-weight elements.
Denote $B_\mu=(R/\PP_\mu)S_\mu^{-1}$.
The algebra $B_\mu$ is generated by $q$-commuting elements
$u_1^{\pm 1},\ldots,u_k^{\pm 1},u_{k+1}^{\pm 1},\ldots, u_{n+m}^{\pm 1} $ with
$$
u_1= x^{(k-1)}_{n-i_1-\cdots-i_k-(k-1)},\ldots,
u_{k-1}= x'_{n-i_1-i_2-1},
u_k= x_{n-i_1}, u_{k+1}=x_{n+1},\ldots, u_{k+m}=x_{n+m}.
$$
For any $s=\{1,2,\cdots,k\}$ we denote
$\psi(s)$ such that 
$$
u_s=x^{(s-1)}_{\psi(s)}.
$$
The generators $u_1,\cdots,u_k,u_{k+1},\ldots,u_{k+m}$ satisfy
the relations 
$$u_iu_j=q'_{ij}u_ju_i\eqno(3.4)$$
and the matrix 
$$\Qb_\mu=(q'_{ij})_{i,j=1}^{k+m}\eqno (3.5)$$
 is a submatrix of $\Qb$.
We denote by $\Sb_\mu$ the corresponding submatrix of $\Sb$.
The algebra
$B_\mu$ admits the filtration
$$B_\mu=B_{\mu,1}\supseteq B_{\mu,2}\supseteq\cdots\supseteq B_{\mu,n}
\supseteq B_{\mu,n+1}\eqno(3.6)$$
with 
$$
B_{\mu,j}= (R_j/({\PP_\mu\bigcap R_j}))S_{\mu j}^{-1},
\qquad
S_{\mu j}=S_\mu\bigcap R_j/(\PP_\mu\bigcap R_j).
$$
Notice that some subalgebras $B_{\mu,j}$ and $B_{\mu,j+1}$
can coincide.
If $B_{\mu,j}\ne B_{\mu,j+1}$, then
$j=\psi(s)$ for some $s\in\{1,2,\ldots,k\}$ and $B_{\mu,j}$ is generated by
$u_s$ and $B_{\mu,j+1}$.
Recall that $R_j$ is a skew extension of $R_{j+1}$ and $x_j$
(see Definition 3.1).
We retain the notation $x_j$ for image of $x_j$ in $B_{\mu,j}$.
We get that $B_{\mu,j}$ is generated by $x_j$ and $B_{\mu,j+1}$.
The elements $u_s$ and $x_j$ are related as follows
$$u_s=x_j+r_{j+1}\eqno(3.7)$$
with $r_{j+1}\in B_{\mu,j+1}$.\\
{\bf Remark 3.7}.\\
1) The ideal $\PP_\mu$ is an intersection of prime ideals $\{\QQ\}$ with
$\QQ\bigcap C=0$. 
Hence, any $c\in C$, $c\ne 0$ is a regular element in $R/\PP_\mu$.
It implies that $R/\PP_\mu$ is a free $C$-module.\\
2) The generators of $S_\mu$ $q$-commute.
Each element $s$ of $S_\mu$ is a monomial of $q$-commuting
generators $u_1,\ldots,u_{k+m}$ (see the stratification process).
The elements $s^l$ lie in the center of
$B_\mu\bmod(q-\eps)$. One can consider the new denominator subset
$S_\mu^l=\{s^l: s\in S_\mu\}$ in $R/P_\mu$. Obviously,
the localization over $S_\mu^l$ coincides with the localization over
$S_\mu$. It follows that we may choose $S_\mu$ such that
it lies in the center of $B_\mu\bmod(q-\eps)$.

Now we are returning to ideals in $R$. 
Let $R$ be a quantum solvable algebra obeying Conditions 3.2-3.4.
Our next goal is to study stratification of prime ideals of specialisations
of $R$.
For any $\lambda\in K$ we consider the specialisation
$\rho_\lambda: R\to R_\lambda=R\bmod(q-\lambda)$.
For any prime ideal $I$ in $R_\lambda$ we denote
$\II=\rho_\lambda^{-1}(I)$. The above ideal $\II$ is prime
in $R$ and $\II\bigcap C=C(q-\lambda)$.
Denote as usual
$\PP_{\mu\lambda}=\PP_\mu\bmod(q-\lambda)=\PP_\mu+R(q-\lambda)/{R(q-\lambda)}$
and $S_{\mu\lambda}$ the image of $S_\mu$ in $R_\lambda/\PP_{\mu\lambda}=
(R/\PP_\mu)\bmod(q-\lambda)$.
Since $\II\bigcap C\ne 0$, we can't apply Theorem 3.6 directly.
The next Example shows that above the stratification 
is not valid for all prime ideals.\\ 
{\bf Example 3.8}. Let $R={\Bbb A}_{q,1}$ be the Quantum Weyl algebra
generated by $x,y$ with the relation $xy-qyx=1$.
Any prime ideal of $R$ with zero intersection with $K[q,q^{-1}]$ 
admits localization over $y$.
The localization of $R$ over $y$ is generated by
$y^{\pm 1}$ and $u=(q-1)yx+1$ related $yu=q^{-1}uy$.
Consider specialisation $q\mapsto\eps=-1$.
The algebra $R_{-1}$ is generated by $x,y$ related
$xy+yx=1$.
Let $I$ be the kernel of irreducible representation
$$
x\mapsto\left(\begin{array}{cc}0&0\\1&0\end{array}\right)
\qquad 
y\mapsto\left(\begin{array}{cc}0&1\\0&0\end{array}\right)
  $$
The ideal $I$ of $R_{-1}$ is prime, 
but does not admit localization over $y$. $\Box$\\
{\bf Definition 3.9}. Let $\lambda\in \overline{K}$ be not a root of unity.
We say that $\lambda$ is
is a point of good reduction of the stratification $\Mb=\{\PP_\mu\}$, 
if the following properties hold. \\
Property 1. Every 
ideal $\PP_{\mu\lambda} $ is a semiprime ideal in $R_\lambda$,
$S_{\mu\lambda}$ is a denominator set of $R_\lambda/\PP_{\mu\lambda}$.\\
Property 2. For any prime ideal $I$ in $R_\lambda$
there exists standard ideal $\PP_{\mu\lambda}$ of $R_\lambda$
such that $I\supseteq\PP_{\mu\lambda}$ and 
$S_{\mu\lambda}\bigcap I/\PP_{\mu\lambda}=\emptyset$.
(i.e $\II\supseteq \PP_\mu$ and $S_\mu\bigcap \II/\PP_\mu=\emptyset$).\\
We say that $\lambda$ 
is a point of bad reduction if the statement of
Definition 3.9  are not true for $\lambda$.

For a point of good reduction we denote
$$B_{\mu\lambda}=(R_\lambda/\PP_{\mu\lambda})S_{\mu\lambda}^{-1}
= B_\mu\bmod(q-\lambda)$$ 
and we can study any prime ideal $I$ in $R_\lambda$ as an ideal
of $B_{\mu\lambda}$ for some $\mu$.
 That is we can reduce the study of 
$I$ to the study of an ideal in an algebra of twisted Laurent polynomials.

Let $\eps\in \Lambda$ (see Condition 3.3). That is
 $\eps$ is a primitive $l^{\rm{th}}$ root of unity
and $x_1^l,\ldots,x_n^l$ lie in the center $Z_\eps$ of $R_\eps$.
Denote by $\DD_0$ the set of derivations
$\DD_{x_1^l},\ldots,\DD_{x_{n+m}^l}$ defined in Section 2.
We say that an ideal $I$ of $R_\eps$ 
is $\DD_0$-stable if $I$ is stable with
respect to all derivations of the set $\DD_0$. \\
{\bf Definition 3.10}. As above $\lambda=\eps$ is 
a primitive $l^{\rm{th}}$ root of unity and 
$x_1^l,\ldots,x_{n+m}^l$ lie in the center of $R_\eps$.
 We say that $\eps$ 
is a point of good reduction of stratification $\Mb=\{\PP_\mu\}$, 
if 
Property 1 of above Definition 3.9 holds and Property 2 holds for any prime
$\DD_0$-stable ideal $I$ of $R_\eps$. \\
We say that $\eps\in\Lambda$ 
is a point of bad reduction if the statement of
Definition 3.10 is not true for $\eps$.
We say further that the property ${\cal E}(\lambda)$ 
is true for almost all $\lambda\in\overline{K}$ 
if there exists a finite set $F$ such 
that this property is true for all elements $\lambda\in\overline{K}-F$. \\  
{\bf Lemma 3.11.}  Let $B'$ be a free $C$-module and an unitary 
$C$-algebra generated by
$q$-commuting elements $u_1,\ldots,u_M$ and $\II$
is a prime ideal in $B'$.
Then the image of 
each $u_i$ in $B'/\II$ is either zero
or a regular element.\\
{\bf Proof}.
If the generator $u_{i}$ is a zero divisor in $B'/\II$, 
then there exists $b\notin I$ such that $u_{i}b\in\II$.
Then for all $j$ we have 
$u_{i}u_{j}b=q_{ij}u_{j}u_{i}b\in\II$.
Hence, $u_{i}B'b\in\II$. Since $\II$ is prime, $u_{i}\in \II$.
$\Box$\\
{\bf Lemma 3.12}. 
Let $R$ be a free $C$-module and an unitary $C$-algebra.
Let $\PP$ be a semiprime ideal in $R$ with $\PP\bigcap C=0$
and $\PP_\eps=(\PP+(q-\eps)R)\bmod(q-\eps)$. 
We claim the following.\\
1) $\PP_\eps$ is  $\DD$-stable;\\
2) Let $u_\eps=u\bmod(q-\eps)\in Z_\eps$ and 
$\DD_u$ the quantum adjoint action of $u$.
We denote by $\bar{u}$ the image of $u$ in $B=R/\PP$. As usual 
$B_\eps=B\bmod(q-\eps)$. We denote by $\overline{\DD}_{\bar{u}}$ 
the quantum adjoint action of $\bar{u}$ on $B_\eps$ 
via reduction $\bmod(q-\eps)$.
Then 
$ \DD_u :B_\eps\to B_\eps$ 
coincides with $\overline{\DD}_{\bar{u}}$;\\
3) 
Denote by $Z'$ the image of the center $Z_\eps$ of $R_\eps$ 
under the map
$R_\eps\to R_\eps/{\PP_\eps}=B_\eps$.
The subalgebra $Z'$ lies in the center $Z(B_\eps)$. We equip $Z'$ with the 
Poisson structure
 inherited by the Poisson structure in $Z_\eps$.
 We consider the Poisson structure of 
 $Z(B_\eps)$ defined via reduction
$B_\eps=B\bmod(q-\eps)$ (Section 2).
 We assert that $Z'$ is a Poisson subalgebra of $Z(B_\eps)$, i.e.
the Poisson bracket of two elements in $Z'$ 
 coincides with 
 their Poisson bracket
as elements $Z(B_\eps)$.\\
 {\bf Proof}.1)
 For any $a\in \PP_\eps=(\PP+R(q-\eps))\bmod(q-\eps)$ one can choose 
$\underline{a}\in \PP$ such that $a=\underline{a}\bmod(q-\eps)$. 
Let $u\in R$ such that $u_\eps=u\bmod(q-\eps)\in Z_\eps$, then
$b=u\underline{a}-\underline{a}u\in R(q-\eps)\bigcap\PP=(q-\eps)\PP$.
Then $b\in (q-\eps)\PP$ and $\DD_u(a)\in \PP_\eps$.\\
The claim 2) is obtained by reduction $u\ua-\ua u\in (q-\eps)R$
(with $\ua\in R$) modulo $\PP$.
The claim 3) is a  corollary of claim 2). $\Box$\\
{\bf Notation 3.13}.
Let $\eps\in\Lambda$ and $R$ obey Conditions 3.2-4.
We denote
$$\theta_i = \frac{\tau_i^l-\mbox{id}}{q-\eps} :R_\eps\to R_\eps$$
for $0\le i\le n+m$ and $\Theta=\rm{span}\{\theta_1,\ldots,\theta_{n+m}\}$.
For any ideal $I$ of $R_\eps$ we denote by $I_\Theta$ the greatest 
$\Theta$-stable ideal in $I$.

Notice that any $\theta_i$ is a derivation. Whence 
$\Theta:Z_\eps\to Z_\eps$.  
We shall refer
a common eigenvector for all $\theta_1,\ldots,\theta_{n+m}$ as a 
$\Theta$-weight vector.\\
{\bf Proposition 3.14}. Let $\eps\in\Lambda$ and $R$ obey Conditions 3.2-4. \\
1) Let $\tau\in\{\tau_1,\ldots,\tau_{n+m}\}$ and let $\theta$ be
the correspoding element of $\{\theta_1,\ldots,\theta_{n+m}\}$.
Let $u\in R$ such that $u_\eps=u\bmod(q-\eps)$ and $\tau(u)=q^mu$. 
Then \\
a) $\theta(u_\eps)=\um u_\eps$,\\
b) $\tau\DD_u=\eps^m\DD_u\tau$,\\
c) $\theta\DD_u=\DD_u(\theta+\um)$;\\
2) if an ideal $I$ of $R_\eps$ is $\DD$-stable (cor.$\DD_0$-stable), 
then $I_\Theta$ is also $\DD$-stable (cor.$\DD_0$-stable);\\
3) let $I$ be an ideal in $R_\eps$ , $\rho_\eps:R\to R_\eps$ be the
specialisation and $\II=\rho^{-1}_\eps(I)$. Then 
$\rho^{-1}_\eps(I_\Theta)=\II_H$ where $\II_H$ is the 
greatest $H$-stable ideal in $\II$.
In particular, 
if an ideal $I$ of is $\Theta$-stable, then the ideal $\II$
is $H$-stable.\\
{\bf Proof}.\\
1) The proof of a) easy (see Section 2).
For any  $a\in R_\eps$ choose $\ua\in R$ such that $a=\ua\bmod(q-\eps)$.
We have
$$
\tau\left(\frac{u\ua-\ua u}{q-\eps}\right) =
q^{m}\left(\frac{u\tau(\ua)-\tau(\ua)u}{q-\eps}\right).$$
This yields $\tau\DD_u(a)=\eps^m\DD_u\tau(a)$ and proves b).\\
By direct calculations,
$$
\left(\frac{\tau^l-\mbox{id}}{q-\eps}\right)
\left(\frac{u\ua-\ua u}{q-\eps}\right)=
q^{ml}\left(\frac{u\tau^l(\ua)-\tau^l(\ua)u}{(q-\eps)^2}\right)-
\frac{u\ua-\ua u}{(q-\eps)^2}=
$$
$$
q^{ml}\left({u
\frac{\tau^l-\mbox{id}}{q-\eps}(\ua)-
\frac{\tau^l-\mbox{id}}{q-\eps}(\ua)u}\right)\cdot\frac{1}{q-\eps}
+
\frac{q^{ml}-1}{q-\eps}\cdot\frac{u\ua-\ua u}{q-\eps}.$$
Putting $q=\eps$, we get 
$\theta\DD_x(a)=\DD_x\theta(a)+\underline{m}\DD_x(a)$. This proves c).\\
2) An element $a$ lies in $I_\Theta$ iff $a\in I$ and for any 
$i_1,\ldots,i_k\in\{1,\ldots,n+m\}$ the element
$\theta_{i_1}\cdots\theta_{i_k}(a)$ also lies in $I$. 
Suppose that $I$ is $\DD$-stable. The center $Z_\eps$ is spaned
 by $\Theta$-weight elements. Let $u_\eps$ be a $\Theta$-weight element
  in $Z_\eps$. Choose a $H$-weight element $u\in R$ such that 
  $u_\eps=u \bmod(q-\eps)$. 
  Then for any $a\in I_\theta$ we have
 $$
 \theta_{i_1}\cdots\theta_{i_k}(\DD_u a)=
 \DD_u(\theta_{i_1}+\um_{i_1})\cdots(\theta_{i_k}+\um_{i_k})(a)\subset 
 \DD_u(I)\subset I.
 $$
 It follows $\DD_u(a)\in I_\Theta$. Whence $I_\Theta$ is $\DD$-stable.
 Similar for $\DD_0$-stable ideals.\\
3) For any ideal $I$ of $R_\eps$
the ideal $I_\Theta$ is a span of all $\Theta$-weight elements in $I$.
Similar 
the ideal $\II_H$ is a span of all $H$-weight elements in 
$\II=\rho^{-1}_\eps(I)$.
For any $\Theta$-weight element $a_\eps\in I$ one can choose  
a $H$-stable element $a\in \II$ such that 
$a_\eps= \rho_\eps(a)=a \bmod(q-\eps)$. 
We have $\rho^{-1}_\eps(I_\Theta)\subseteq\II_H$.
On the other hand, if an element $a$ is $H$-weight, then
the element $a_\eps$ is $\Theta$-weight.
This proves the inverse containment.$\Box$\\\   
{\bf Theorem 3.15}.\\
 1) There are at most finite number of points of bad reduction among
non roots of unity.\\
2) Denote by $\Lambda_1=\{\eps\in\Lambda: l{\mbox{~is~relatively~ 
prime~with~}}
s_1,\ldots, s_n\}$ (see Condition 3.2).
There are at most finite number of points of bad reduction in $\Lambda_1$.\\
{\bf Proof}.\\ 
{\bf Step 1}. Suppose that $\lambda$ is not a root of 1.
First consider the semiprime $H$-stable (see Conditon 3.3) ideal
$\PP_{i_1}$ defined in (3.2) as an intersection of
prime ideals $\QQ\in X_1'$.
These prime ideals are $H$-stable [P2, Prop.2.1].
For almost all $\lambda$ each ideal $\QQ_\lambda$
is semiprime [P2, Prop.2.7]. Recall that $x_{n-i_1}\notin \QQ$ for
$\QQ\in X_1'$. For almost all $\lambda$  the element $x_{n-i_1,\lambda}$ 
is regular in $R_\lambda/\QQ_\lambda$ [P2,Prop.2.7]. 
This proves that $\PP_{i_1,\lambda}$ is semiprime
and the element $x_{n-i_1,\lambda}$ is regular
in $R_\lambda/\QQ_\lambda$ for each $\QQ\in X_1'$.
The set $ S_{1*}=\{x_{n-i_1}^m: m\in\Nb\}$
obeys the Ore condition in $R/\PP_{i_1}$ [P2,Theorem 3.4].
Then $S_{1*,\lambda}$ obeys Ore condition in $R_\lambda/\PP_{i_1,\lambda}$.
We conclude that $S_{1*,\lambda}$ is a denominator set
and consider the localization $R^{(1)}_\lambda$ of $R_\lambda/{\PP_{i_1,\lambda}}$
over $S_{1*,\lambda}$ .

For almost $\lambda$ the elements of denomimator set $N_1\subset C$
(see the stratification process) is non zero modulo $q-\lambda$.
We can consider localization over $N_1$ modulo $q-\lambda$.
Following the stratification process, we see that all ideals
$\PP_{\mu\lambda}$ are semiprime and $S_{\mu\lambda}$ is a denominator
set in the factor algebra of $R_\lambda$ over $\PP_{\mu\lambda}$. This proves
Property 1 of Definition 3.9.

Notice that, since $\lambda$ is not a root of unity, then 
all prime ideals in $R$ are completely prime 
[GL,Theorem 2.3].

Consider the nonempty Zariski-open set $\Ob_{st}\subset \overline{K}$ defined as follows
(i.e. $ \overline{K}-\Ob_{st}$ is finite).
A point $\lambda$ lies in $\overline{K}-\Ob_{st}$ if $\lambda$ 
annihilate some of ideal in $ X''_m\bigcap C$ (see the stratification process)
or some element in 
$S_\mu\bigcap C$.
The Property 2 of Definition 3.9 is true for all $\lambda$ in $\Ob_{st}$
[P2,Theorem 3.10]. \\
{\bf Step 2}. Suppose that $\lambda=\eps\in \Lambda_1$. 
We assume that $\lambda\in \Ob_{st}$.
The proof of Property 1) of Definition 3.10 is similar to Step 1.
We are going to prove Property 2) of the above definition.
Let $I$ be a prime $\DD_0$-stable ideal in $R_\eps$.
As usual $\II=\rho_\eps^{-1}(I)$  with $\II\bigcap C=(q-\eps)C$.
Denote $\JJ=\II_H$ and $J=I_\Theta$. 
By Proposition 3.14, $J$ is a $\DD_0$-stable ideal and 
$\JJ=\rho_\eps^{-1}(J)$. By 
[MC-R,14.2.3] and [D,3.3.2], the ideal $J$ is prime.  

Recall that 
Theorem 2.12 asserts that for all $j$ the ideal $J_j=J\bigcap R_{j\eps}$
is prime in $R_{j\eps}$. It follows that $\JJ_{j}=\JJ\bigcap R_j$ 
is a prime
ideal in $R_j$. 

Suppose that some subset $L_1=\{x_n,\ldots,x_{n-i_1-1}\}\subset \JJ$
and $x_{n-i_1}\notin \JJ$. It follows that $L_1\subset \II$.
Since all generators, in particular $x_{n-i_1}$,
 are $H$-weight elements, then $x_{n-i_1}\notin \II$. 
The ideal $\JJ$
contains a minimal prime ideal, say $\QQ$, in the set $X_1$   
(see the stratification process).
Since $x_{n-i_1}\notin \JJ$, then $\QQ\in X_1'\bigcup X_1''$.
If $\QQ\in X_1''$, then $\eps$ annihilates $Q\bigcap C$.
This contradicts the claim $\eps\in\Ob_{st}$.
We get $\QQ\in X_1'$ and 
 $\JJ$ contains the $H$-stable ideal $\PP_{i_1}$. 
We saw above that $\JJ_{n-i_1}=\JJ\bigcap R_{n-i_1}$ 
is a prime ideal in $R_{n-i_1}$. The algebra 
$$
R'_{n-i_1}=\frac{R_{n-i_1}}{\PP_{i_1}\bigcap R_{n-i_1}}
$$ is generated by 
$q$-commuting elements. The ideal 
$$
\JJ'_{n-i_1}=\frac{\JJ_{n-i_1}}{\PP_{i_1}\bigcap R_{n-i_1}}$$ is a prime
ideal in $R'_{n-i_1}$. 

We apply Lemma 3.11 for $\JJ'_{n-i_1}$.
Since $x_{n-i_1}\notin \JJ$, the image of this element is regular in
$R'_{n-i_1}/\JJ'_{n-i_1}$. Therefore, $\JJ'_{n-i_1}$
has empty intersection with
denominator set $S_{1*}=\{x^m_{n-i_1}\}_{m\in\Nb}\subset R/\PP_{i_1}$.
By [D,3.6.17] or [MC-R,2.1.16], the ideal $\JJ'$ admits localization 
over $x_{n-i_1}$. Since $S_{1*}$ consists of $H$-weight elements,
the ideal $\II'=\II/\PP_{i_1}$ also has empty intersection with $S_{1*}$ and
admits localization on $S_{1*}$. 
 
 Since $\eps\in \Ob_{st}$, then $\eps$ does not annihilate any
 element of $S\bigcap C$. Hence $\eps$ does not annihilate $N_1$
 (see the stratification process).
 The ideals $\JJ'$ and $\II$ admit localization over $N_1$ (see 
stratification process).

Suppose that  $\JJ'$ contains some subset 
$L_2=\{x'_{n-i_1-1},\ldots,x'_{n-i_1-i_2}\}$ and does not contain
$x'_{n-i_1-i_2-1}$. It follows that $L_2\subset \II'$.
Since the generators are $H$-weight elements,
 $x'_{n-i_1-i_2-1}\notin\II'$.
As above the claim  $\eps\in\Ob_{st}$ implies $\JJ\supseteq \PP_{i_1i_2}$ 
(see also [P2,Theorem 3.10]).
As above the ideal
$\JJ''=\JJ/ \PP_{i_1i_2}$ has empty intersection with 
denominator set $\{(x'_{n-i_1-i_2-1})^m\}_{m\in\Nb}$.
The ideal $\JJ''$ admits localization over $x'_{n-i_1-i_2-1}$.
The ideal $\JJ''/\PP_{i_1}S^{-1}_{2*}$ admits 
localization over $N_2$. The same is true for 
$\II''=\II/ \PP_{i_1i_2}$.
And so on.

Finally, we prove that ideal $\II$ contains a unique ideal 
$\PP_\mu$ such that
$S_\mu\bigcap \II/{\PP_\mu}=\emptyset$.
Then $I\supseteq \PP_{\mu\eps}$ and  $S_{\mu\eps}\bigcap I/{\PP_{\mu\eps}}
=\emptyset$. Notice that, by Lemma 3.12, the ideal 
$\PP_{\mu\eps}$ is $\DD$-stable and, by [P2, Prop.2.7], semiprime for
almost all $\eps$. 
$\Box$

Since $B_\mu$ is generated by $q$-commuting elements, we shall 
recall some well known statements on algebras
of twisted Laurent polynomials.
Let $\Qb$ and $\Sb$ be as above
$M\times M$-matrices with $q_{ij}=q^{s_{ij}}$ .
Let $B=B_{\Qb}$ be a free $C$-module and a unitary $C$-algebra generated by 
$u_1^{\pm 1},\ldots,u_M^{\pm 1}$ with relations $u_iu_j=q_{ij}u_ju_i$.
The algebra $B$ is a factor algebra of an algebra of twisted Laurent
plynomials.
One can decompose $B=A\otimes Z(B)$ where $A$ is an algebra of twisted
Laurent polynomials and $Z(B)$ is the center of $B$.
The generators $y^{\pm 1}_1,\ldots,y^{\pm 1}_{2r}$ of $A$ are monomials of
 $u_1^{\pm 1},\ldots,u_M^{\pm 1}$ and obey the relations
$$
y_1y_2=q^{d_1}y_2y_1,\ldots, 
y_{2r-1}y_{2r}=q^{d_r}y_{2r}y_{2r-1}\eqno (3.8)
$$ 
with positive integers 
$d_1,\ldots,d_{2r}$. 
All other pairs of generators commute.
Here $2r$ is the rank $\Sb$. We shall refer $2r$ as the rank of $B$.  
The center $Z(B)$ is generated by monomials.
Any ideal of $B$ is generated by its intersection with the center $Z(B)$.
Any prime ideal of $B$ is completely prime.

One can consider the specialisation $B_\eps$ of $B$ at primitive
$l^{\rm{th}}$ root of unity $\eps$.
We obtain $B_\eps=A_\eps\otimes Z(B)_\eps$.\\
{\bf Lemma 3.16}. 
Suppose that $l$ is relatively prime with the elementary 
divisors of matrix $\Sb$.
Then 
1) the center $Z(B_\eps)$ is generated as $K$-algebra 
by $u_{1\eps}^{\pm l},\ldots, u_{M\eps}^{\pm l}$ and the reduction $Z(B)_\eps$
of $Z(B)$ modulo $q-\eps$;\\
2) if the field $K$ is algebraically closed, then the dimension of
every irreducible representation of $B_\eps$ equals to 
$l^r$ with $r=\frac{1}{2}\mbox{rank}(\Sb)$.\\
{\bf Proof}. Consider another system 
$y^{\pm 1}_1,\ldots,y^{\pm 1}_{2r},z^{\pm 1}_1,\ldots,z^{\pm 1}_t$, $2r+t=M$ 
of generators of $B$ with the relations $y_1y_2=q^{d_1}y_2y_1,\ldots, 
y_{2r-1}y_{2r}=q^{d_r}y_{2r}y_{2r-1}$. All other pairs of generators commute.
Since $(l,d_i)=1$ for all $i$, then the $l^{\rm{th}}$ powers of  
generators and $z_1,\ldots,z_t$ 
generate the center of $B_{\eps}$. This proves 1).

Any irreducible representation $\pi$ of $B_\eps$ determines the central 
character $\chi:Z(B_\eps)\to K$:
 $$\chi(y_i^l)=\gamma_i\in K^*,\qquad
\chi(z_j)=\alpha_j\in K^*.$$
The factor algebra $B_\eps/m(\chi)B_\eps$ is isomorphic to a
matrix algebra over $K$.
This proves that any irreducible representation uniquely determines 
by its central 
character
and has the dimension of statement 2).
$\Box$

Further we denote by $\overline{\DD}_u$ the quantum adjoint action of
an element $u\in B$ with $u_\eps\in Z(B_\eps)$.
Clearly, $\oDD_u(Z(B)_\eps)=0$.
For $u_i^l$ we calculate 
$$\overline{\DD}_{u_i^l}(u_j)=s_{ij}l\eps^{-1}u_ju_i^l,$$
$$ \{u_{i\eps}^l,u_{j\eps}^l\}=s_{ij}l\eps^{-1}u_{i\eps}^lu_{j\eps}^l.$$
We say that an ideal of $B_\eps$ is $\oDD$-stable 
if it is stable with respect to all derivations $\DD_u$ with 
$u_\eps\in Z(B_\eps)$. We say that an ideal of $B_\eps$ 
is $\oDD_0$-stable if it is
stable with respect to all derivations $\DD_{u_i^l}$ with $1\le i\le M$.\\
{\bf Lemma 3.17}. Let $l$ be as in Lemma 3.16.
1) Any $\oDD_0$-stable ideal of $B_\eps$ is generated by its 
intersection with $Z(B)_\eps$;
2) Any Poisson ideal in $Z(B_\eps)$ is generated by its intersection with
$Z(B)_\eps$;
3) For $K=\Cb$, each symplectic leaf of ${\rm{Maxspec}}Z(B_\eps)$
 has dimension $2r$, it is closed in Zariski topology and is defined by
 the ideal $\sum_j^t Z(B_\eps)(z_j-\alpha_j)$ where $z_1,\ldots,z_t$ generate
 $Z(B)_\eps$ and $\alpha_1,\ldots,\alpha_t$ are some complex numbers.   \\
{\bf Proof}. Choose the above generators $y_1,\ldots,y_{2r}$ in $B$ over
$Z(B)$. Recall that each $y_i$ is a monomial of $u_1,\ldots,u_M$.
Consider the derivations $\Phi_i=y_{i\eps}\oDD_{y_i^l}$ of $B_\eps$.
Let $\Phi$ be the linear subspace spanned by $\Phi_i$, $1\le i\le 2r$.
If $I$ is $\oDD_0$-stable ideal, then it is $\Phi_i$-stable for all 
$1\le i\le 2r$ (call $\Phi$-stable).
The monomials $Y_{\overline{m}}
=y_1^{m_1}\cdots y_{2r}^{m_{2r}}$ with $\overline{m}=(m_1,\ldots, m_{2r})$ 
form
a $Z(B)$-basis of $B$.
The images $Y_{\overline{m}\eps}$ of $Y_{\overline{m}}$ in
$B_\eps$ form a $Z(B)_\eps$-basis of $B_\eps$. 
All $Y_{\overline{m}\eps}$ are $\Phi$-eigenvectors:
$$\Phi_i(Y_{\overline{m}\eps})=
\pm d_{i\pm 1}m_il\eps^{-1}Y_{\overline{m}\eps}$$
with different $\Phi$-weights. This implies that
any $\Phi$-stable subspace is generated  by elements 
$Y_{\overline{m}\eps}z$ with $z\in Z(B)_\eps$.
Since the elements $Y_{\overline{m}\eps}$ are invertible, any $\Phi$-stable
ideal is generated  by its intersection with $Z(B)_\eps$. Whence the
same is true for any $\oDD_0$-stable ideal in $B_\eps$. This proves 1). 
It implies 2) and 3). 
$\Box$\\
{\bf Corollary 3.18}. Let $l$ and $B$ be as above.
Any prime
$\oDD_0$-stable (in particular, $\oDD$-stable) ideal $I$ of $B_\eps$ is 
completely prime. \\
{\bf Proof}. 
Let $I$ be a prime $\oDD_0$-stable ideal in $B_\eps$.
By the above Lemma, $I$ is generated by its intersection with $Z(B)_\eps$. 
The intersection
$I\bigcap Z(B)_\eps$ is a completely prime ideal in $Z(B)_\eps$ and
$$B_\eps/I=A_{\eps}\otimes \frac{Z(B)_\eps}{I\bigcap Z(B)_\eps}.\eqno(3.9)$$
This proves the claim. $\Box$

 Let $B=B_\Qb$  be as above. 
Consider the filtration
of $B$ by subalgebras $B_j$ generated by $u^{\pm 1}_j,\ldots,u^{\pm 1}_M$:
$$ B=B_1\supseteq\cdots\supseteq B_j\supseteq\cdots \supseteq B_M
\supseteq B_{M+1}=C.$$
The center $Z(B)$ has a filtration 
$$ Z(B)=Z(B)_1\supseteq\cdots\supseteq Z(B)_j\supseteq\cdots
\supseteq Z(B)_M\supseteq
Z(B_{M+1})=C$$
where
$Z(B)\bigcap B_j= Z(B)_j$.
Similarly, $B_\eps=B\bmod(q-\eps)$ has the filtration $B_{j\eps}$    
 with $Z(B_\eps)\bigcap B_{j\eps}= Z(B_\eps)_j$. Clearly, 
 $Z(B)_{j\eps}\subset Z(B_\eps)_j$.\\
{\bf Lemma 3.19}. Let $\Sb_j$ be a 
$S$-submatrix of $(s_{i\beta})$, $1\le i\le M$, $j\le\beta\le M$.
Suppose that $l$ is relatively prime with elementary
divisors of $\Sb_j$ .
If  $a=u_j^{n_1}\cdots u_M^{n_M}$, $ n_i\in\Zb$ lies in the center
$Z(B_{\Qb,\eps})$
then there exist the $m_\beta\in\Zb$, $j\le\beta\le M$, 
$m_\beta=n_\beta\bmod(l)$ such that 
$a=u_j^{m_1}\cdots u_M^{m_M}$ lies in the center
$Z(B)$.\\
{\bf Proof}. For any $i$ and any monomial $a=u_j^{n_1}\cdots u_M^{n_M}$ 
we have 
$$u_ia=u_i(u_j^{n_1}\cdots u_M^{n_M})=
q^{\sum s_{i\beta}n_\beta} u_j^{n_1}\cdots u_M^{n_M}u_i=
q^{\sum s_{i\beta}n_\beta}au_i.
$$
  A monomial $a$ lies in $Z(B_{\eps})$ if
$(n_j,\ldots,n_M)$ is a solution of the system of equations

$$\sum_{\beta=j}^Ms_{i\beta}n_\beta=0\bmod(l),\qquad 1\le i\le M.$$
Under the requirements on $l$ any solution modulo $l$ is a reduction modulo
$l$ of a solution over $\Zb$.$\Box$\\
{\bf Definition  3.20}. 
We say that a primitive $l^{\rm{th}}$ root of unity $\eps$ is admissible 
if 
the elements $x_1^l,\ldots,x_{n+m}^l$ lie in the center of $R_\eps$
(i.e. $\eps\in\Lambda$),
$\eps$ is a 
point of good reduction
(see Definition 3.10), and 
$l$ is 
relatively prime with all minors of integer matrix
$\Sb$ (Definition 3.1) and with all integers $s_1,\ldots,s_n$ (Condition 3.2). 

 Suppose that $\eps$ is admissible.
Consider the reduction of ideals $\PP_\mu$ and denominator sets
$S_\mu\subset R/{\PP_\mu}$ modulo $q-\eps$.
We denote $B=B_\mu=(R/{\PP_\mu})S_\mu^{-1}$ and $B_\eps=B\bmod(q-\eps)$.
We get  $\DD$-stable semiprime ideals $\PP_{\mu\eps}$ in $R_\eps$. 
Recall that $B_\eps$ is generated by $q$-commuting elements
$u_1^{\pm 1},\ldots ,u_{k+m}^{\pm 1}$ (see (3.4),(3.5)).

One can extend the $\DD$-action from $R_\eps$ to $B_\eps$.
The algebra $B_\eps$ admits also the quantum adjoint action $\overline{\DD}$ 
of its center $Z(B_\eps)$.
Our next goal is to compair these
 two actions.
 We consider $Z(B_\eps)\overline{\DD}$ the $Z(B_\eps)$-module spanned by 
$z_\eps\overline{\DD}_b$ with $b\in B$, $b_\eps=b\bmod(q-\eps)\in Z(B_\eps)$
and $z_\eps\in Z(B_\eps)$.
 We consider the $Z(B_\eps)$-submodule
$\mbox{ad}_{B_\eps}$ of inner derivations in $B_\eps$.
We denote by $Z(B_\eps)\DD_0$ the $Z(B_\eps)$-submodule
$$Z(B_\eps)\DD_{x_1^l}+\cdots+Z(B_\eps)\DD_{x_n^l}$$
of derivations of $B_\eps$.\\
{\bf Proposition 3.21}. As above $R$ is a quantum solvable algebra obeying 
Conditions 3.2-3.4 and $\eps$ is admissible.
Then
$Z(B_\eps)\DD_0 = Z(B_\eps)\overline{\DD}\bmod(\mbox{ad}_{B_\eps})$.\\
{\bf Proof}. 
First $Z(B_\eps)\DD_0 \subseteq Z(B_\eps)\overline{\DD}$
(by Lemma 3.12, the derivation $\DD_{x_i^l}:B_\eps\to B_\eps$ coincides with 
$\overline{\DD}$ derivation of the image of $x_i^l$ in $B$).
To conclude the proof we are going to prove 
$Z(B_\eps)\DD_0$ contains $\overline{\DD}$ modulo inner derivations.

 Consider the filtration of $B=B_\mu$ (see (3.6)).
 Recall that by $Z(B_\eps)_j$ we denote the intersection of $Z(B_\eps)$ with
$B_{j\eps}$.
We consider 
$$\overline{\DD}_j= \overline{\DD}_{Z(B_\eps)_j}=
\{\overline{\DD}_u: u\in B, u_\eps=u\bmod(q-\eps)\in Z(B_\eps)_j\}.$$
To within an inner derivation, we may
choose above $u\in B_j$.
We denote
$Z(B_\eps)_j\DD_{0j}= Z(B_\eps)_j\DD_{x_j^l}+\cdots+Z(B_\eps)_j\DD_{x_n^l}$.
We 
shall prove by induction on $j$ that $\overline{\DD}_j$ equals to 
$Z(B_\eps)_j\DD_{0j}$ modulo inner derivations.

The statement is obviously true for $j=n+1$.
Assume that the statement is true for 
$j+1$. Our goal is to prove it for $j$.   

Let $B_j\ne B_{j+1}$.
The algebra $B_j$ is generated by $B_{j+1},x_j$
and also by $B_{j+1},u_s$. The elements $u_s$ and $x_j$ are related
as (3.7). To simplify the notations we put $u_s=u$.
Recall that $u_\eps^l\in Z(B_\eps)$ and $x_{j\eps}^l\in Z(R_\eps)$. \\
1) Suppose that there exists a monomial 
$z=u^mb\in Z(B)_j$ with $m\ne 0$ and $b\in B_{j+1}$.
We denote by $m_0$ the least positive integer 
with the following property.
The element $v=u^{m_0}b_0$ modulo $q-\eps$ lies in $Z(B_{\eps})_j$
for some $b_0\in B_{j+1}$. 
Here one may consider $b_0$ to be a monomial.  
The monomial $b_0$ is unique up to multiplicaton by a monomial which modulo
$q-\eps$ lies in
$Z(B_\eps)_{j+1}$.

The algebra $Z(B_{\eps})_{j}$ is generated by $v_\eps=v\bmod(q-\eps)$ and 
$Z(B_{\eps})_{j+1}$.
Hence, $m_0$ divides $m$. Put $m=pm_0$.
We obtain $z=v^pa$ with a monomial $a\in B_{j+1}$ and $p\ne 0$.
Notice that  $a\bmod(q-\eps)$ lies in $Z(B_\eps)_{j+1}$.

Since $z\in Z(B)_j$, we have $\overline{\DD}_z=0$.
We obtain
$$0=\overline{\DD}_z=pv_\eps^{p-1}a_\eps\overline{\DD}_v + 
v_\eps^p\overline{\DD}_a.$$
Then $\overline{\DD}_v=-p^{-1}v_\eps a_\eps^{-1}\overline{\DD}_a$.
 By assumption of induction, the derivation $\overline{\DD}_a$ lies in
 $Z(B_\eps)_{j+1}\DD_{0,j+1}$ modulo inner derivations. 
 This proves the statement.\\
2) Suppose that the condition $z=u^mb\in Z(B)_j$ implies
$m=0$. This proves that
$Z(B)_j= Z(B)_{j+1}$.
Consider the element $v=u^{m_0}b_0$ from point 1). The algebra $Z(B)_\eps$ is 
generated by $v_\eps$ and $Z(B_\eps)_{j+1}$.
By the choise, $m_0\le l$. If $m_0<l$, then by Lemma 3.19  
there exists a monomial $u^mb\in Z(B)_j$ with $m_0=m\bmod(l)$.
By the assumption of point 2), $m=0$. Whence, $m_0=l$ and
$Z(B_\eps)_j$ is generated by $u_\eps^l$ and $Z(B_\eps)_{j+1}$.
By (3.7), $u=x_j+r_{j+1}$ with $r_{j+1}\in B_{j+1}$.
Then $ u^l=x_j^l+F(u)$.
Where $F(u)$ a polynomial of $u$ of degree $\le{l-1}$ and $F(u)_\eps$ lies in
$Z(B_\eps)_{j}$. The algebra $Z(B_\eps)_{j}$ is generated by monomials.
 Applying Lemma 3.19, we deduce 
 that $F(u)\bmod(q-\eps)\in Z(B_\eps)_{j+1}$. That is $F(u)=(q-\eps)r+c$
 for some $r\in B_j$ and $c\in B_{j+1}$. The element $c$ 
 modulo $q-\eps$ lies in 
 $Z(B_\eps)_{j+1}$.   
We obtain
$$\overline{\DD}_{u^l}=\overline{\DD}_{x_j^l} + \overline{\DD}_{F(u)}=
\DD_{x_j^l}+\ad_{r_\eps}+\oDD_c \in
\DD_{x_j^l}+ \ad_{r_\eps}+ \overline{\DD}_{Z(B_\eps)_{j+1}}.$$
By the induction hypothesis, we prove the claim. 
 $\Box$.\\
{\bf Corollary 3.22}. Let $L$ be either a  $\DD_0$-stable
ideal in $B_\eps$ or  $\DD_0$-stable ideal in $Z(B_\eps)$, then 
$L$ is  $\overline{\DD}$-stable.\\  
{\bf Notations 3.23}.\\
1) Let $R$ be a domain and a $C$-algebra.
Let $I$ be an ideal in $R_\eps=R\bmod(q-\eps)$.
We denote by $I_\DD$ the greatest 
 $\DD$-stable (i.e. stable with respect to quantum adjoint action)
ideal in $I$. \\
2) Let $R$ be quantum solvable algebra such that 
$x_1^l,\ldots,x_n^l$ modulo $q-\eps$ lie in the center of 
$R_\eps$ (i.e. $\eps\in\Lambda$).
We denote by $I_{\DD_0}$ the greatest $\DD_0$-stable ideal in $I$.\\
{\bf Theorem 3.24}. 
As above $R$ is a quantum solvable algebra
obeying Conditions 3.2-3.4, and $\eps$ is admissible.\\
1) 
Any prime  $\DD_0$-stable (in particular,  $\DD$-stable) 
ideal $I$ is completely prime.\\
2) For above $I$ skew
field of fractions $\Fract(R_\eps/I)$ is isomorphic to a skew field of an
algebra of twisted polynomials.\\
{\bf Proof}.
 There exists standard $\PP_\mu$ such that
$I\supseteq \PP_\mu$ and $I/\PP_{\mu\eps}$ admits 
localization over $S_{\mu\eps}$.
The ideal $I'=(I/\PP_{\mu\eps})S_{\mu\eps}^{-1}$ of 
$B_{\mu\eps}$ is prime $\DD_0$-stable.
By Corollary 3.22, $I'$ is 
$\overline{\DD}$-stable.
The algebra $B_{\mu\eps}$
 is a factor algebra of an algebra of twisted
Laurent polynomials.
By Corollary 3.18, the ideal $I'$ is completely prime. Whence
$I$ is completely prime. 
This proves 1). The claim 2) follows from the observation
$\Fract(R_\eps/I)=\Fract(B_{\mu\eps}/I')$ and (3.9).
$\Box$\\
{\bf Corollary 3.25}. For any prime ideal $I$ in $R_\eps$ the 
ideals $I_{\DD_0}$ and $I_{\DD}$ are completely prime.\\
{\bf Proof}. For any prime ideal $I$ of $R_\eps$
the ideals $I_{\DD_0}$ and $I_\DD$ are also prime ideals of $R_\eps$
[MC-R,14.2.3],[D,3.3.2]. The Theorem concludes the proof.
$\Box$

\section{On Poisson algebras}

In this section we recall some notions and prove some statements which we 
shall use in the next section.

Let $\Aa$ be a commutative  affine $\Cb$-algebra with a Poisson  bracket.
We shall say that $\Aa$ is a Poisson algebra
and $X=\Ms(\Aa)$ is a Poisson variety.

We consider the following stratification of $X$ [BG]:
$X=X_0\supset X_1\supset\cdots\supset X_m=\emptyset$ with 
$X_i=(X_{i-1})_{\rm{sing}}$. The corresponding ideals 
$I_0=\sqrt{\{0\}}\subset I_1\subset\cdots\subset I_m=\Aa$ are semiprime for
all $i=0,\cdots,m$. These ideals are Poisson [Pl, Corollary 2.4].
The smooth locus $X_i^0=X_i-X_{i-1}$
of $X_i$ is a complex analytic Poisson variety.
It is a disjoint union of symplectic leaves.
Recall that the symplectic leaf $\Omega_x$ containing the point 
$x\in X$ is the maximal connected complex analytic  variety in X
such that $x\in \Omega_x$ and the Poisson bracket is 
nondegenerate at every point of $\Omega_x$. 

For any system of generators $a_1,\ldots,a_n$ we denote by 
$M(a_1,\ldots,a_n)$ the matrix
$(\{a_i, a_j\})_{i,j=1}^n$. For any point $x\in X$ we consider the 
specialisation
$$M_x(a_1,\ldots,a_n)=
(\{a_i, a_j\}(x))_{i,j=1}^n.$$
The rank of the matrix $M_x(a_1,\ldots,a_n)$ depends of $x$ and does
not depends of the choice of generators [V,2.6].
We denote $$ \rk_x\Aa=\rk_xM(a_1,\ldots,a_n).$$ 
Notice that $\rk_x\Aa=\dim \Omega_x$.\\
{\bf Definition 4.1}.
We shall denote by $\rk\Aa$  and call rank of $\Aa$ the maximum of $\rk_x\Aa$,
 $x\in X$.

Let $\Bb$  be a Poisson subalgebra in $\Aa$
and $Y=\Ms(\Bb)$. Denote by $\phi:X\to Y$ the corresponding Poisson 
morphism of algebraic varieties.
For any system of generators $b_1,\ldots,b_m$ of $\Bb$ there
exists a matrix $T$ with entries in $\Aa$ such that
$M_y(b_1,\ldots,b_m)=T^tM_x(a_1,\ldots,a_n)T$ for $y=\phi(x)$.
It follows that $\rk_y\Bb\le\rk_x \Aa$.\\
{\bf Proposition 4.2}.
 Let $\Aa$ and $\Bb$ as above. Suppose that $\Aa$ is finite over
$\Bb$. Then there exists a Zariski-open nonempty subset $U\subset Y$ such that
for all $y\in Y$ and $x\in\phi^{-1}(y)$ there holds $\rk_x\Aa=\rk_y\Bb$.\\
{\bf Proof}. The radical of an algebra and all its minimal prime ideals
are stable with respect to any its  derivation [D, 3.3.3].
It follows that  the radical of Poisson algebra and all its 
minimal prime ideals are Poisson ideals.
It is sufficient to prove the claim in the case 
$\Aa$ is a domain.   

Let $\Aa$  be a domain and $a_1,\ldots,a_n$ be generators of $\Aa$.
Consider the first generator $a=a_1$.
Since $\Aa$ is finite over $\Bb$ there exists a
polynomial 
$f(t)=c_0t^k+c_1t^{k-1}+\ldots+c_k$ with $c_0,c_1,\ldots,c_k\in \Bb$ and 
$c_0\ne 0$ such that $f(a)=0$.
We may assume that $f(t)$ is irreducible over $\Fract(\Bb)$.
This implies that the discriminant 
$\Disc(f)$ is a nonzero element of $\Bb$.
Denote by $U_1= \{y\in Y: \Disc(f)(y)\ne 0, c_0(y)\ne 0\}$ and
let $x\in \phi^{-1}(U_1)$. We have $\phi(x)=y\in U_1$.

Notice that $a(x)$ is a root of the polynomial
$f_{(y)}(t)= c_0(y)t^k+c_1(y)t^{k-1}+\ldots+c_k(y)$.
Since $\Disc f_{(y)}(y)=\Disc(f)(y)\ne 0$, all roots (in particular, $a(x)$)
of $f_{(y)}(t)$ are simple. That is 
$$f'_{(y)}(a(x))\ne 0.\eqno (4.1)$$

Consider the matrix $M_x(a_1,\ldots,a_n,c_0,c_1,\ldots,c_k)$.
Its rank is equal to $\rk_x \Aa=\rk M_x(a_1,\ldots,a_n)$.
For any element $g$ of the set $\{a_1,\ldots,a_n,c_0,c_1,\ldots,c_k\}$ we have
$0=\{f(a),g\}=f'(a)\{a,g\}+\{c_0,g\}a^{k}+\ldots+\{c_k,g\}$.
Then
$$f'_{(y)}(a)(x)\{a,g\}(x)+\{c_0,g\}(y)a^{k}(x)+\ldots+\{c_k,g\}(y).\eqno(4.2)$$
The first row and the first column of above matrix is a span of other
rows and columns. Therefore, 
$$\rk_x M(a_1,a_2\ldots,a_n,c_0,\ldots,c_k)=
\rk_x M(a_2,\ldots,a_n,c_0,\ldots,c_k).\eqno(4.3)$$
We may continue the process by putting $a=a_2$.
At the end we have
$$\rk_x\Aa=\rk_y M(c_0,\ldots,c_k,\ldots,c^{(n)}_0,\ldots,c^{(n)}_{k^{(n)}})$$
where all $c_i^{(j)}$ lie in $\Bb$.
Then  $\rk_x\Aa\le\rk_y\Bb$. Hence, $\rk_x\Aa=\rk_y\Bb $ for all $y$ in some
open subset $U$ and $x\in\phi^{-1}(y)$.\\
{\bf Corollary 4.3}.
Let $\Aa$ and $\Bb$ be as in Proposition. Then $\rk \Aa=\rk \Bb$.\\
{\bf Proof}. 
Denote by $X^0$ (resp. $Y^0$) the subset of all $x\in X$ (resp. $y\in Y$)
such that $\rk_x\Aa$ (resp.$\rk_y \Bb$) is not maximal.
For any $x$ in the intersection of the open subsets
$X-X^0$, $\phi^{-1}(U)$ (see above Proposition), and $\phi^{-1}(Y-Y^0)$, we
have
$\rk \Aa=\rk_x\Aa=\rk_y\Bb=\rk \Bb.$
$\Box$
 
\section{Representations at roots of unity}
Let $K=\Cb$ be the field of complex numbers.
Let $R$ be a quantum solvable algebra obeying Conditions 3.2-3.4.
Suppose that $\eps $ is admissible $l^{\rm{th}}$ root of unity (see Definition 3.20). 
As above, $Z_\eps$ is the center of $R_\eps$. The algebra $Z_\eps$ is a 
Poisson algebra.
It determines the algebraic Poisson manifold
$\MM=\Ms(Z_\eps)$.
 
Let $\pi$ be an irreducible complex representation 
of a quantum solvable algebra 
$R_\eps$.  
Since $R_\eps$ is finite over its center, then $\pi$ has finite dimension.
Denote by $I(\pi)$ the kernel of the representation $\pi$ in $R_\eps$.
For $I(\pi)$ we consider $I(\pi)_\DD$-the greatest  $\DD$-stable ideal in 
$I(\pi)$. 

Denote by $\chi=\chi(\pi):Z_\eps\to K$ the central
character of $\pi$.  
The intersection   $I(\pi)\bigcap Z_\eps$ is the maximal ideal
$m(\chi)$ for $\chi\in Z_\eps$. 

We denote by $\Omega_\chi$ the symplectic leaf of $\chi$ with respect to 
Poisson bracket on $Z_\eps$. 
We consider the greatest  $\DD$-stable ideal $m(\chi)_\DD$ in
$m(\chi)$.
The ideal $m(\chi)_\DD$ is prime [D,3.3.2].
We denote by $\MM_\chi$ the annihilator set $\mbox{Ann}(m(\chi)_\DD)$.
Remark that $\MM_\chi$ is the Zariski closure of the symplectic leaf 
$\Omega_\chi$. 
We shall prove further (Theorem 5.3) 
that $\Omega_\chi$
is open in its closure $\MM_\chi$.
The algebra of regular functions $\FF=\Cb[\MM_\chi]$ 
coincides with $Z_\eps/m(\chi)_\DD$. The Poisson structure on $Z_\eps$
provides the Poisson structure in $\FF$. 

 Notice that an element $a$ lies in $I(\pi)_\DD$
iff  for any system of elements $u_1,\ldots,u_s\in R$,
obeying $u_{1\eps},\ldots, u_{s\eps}\in Z_\eps$, 
the element
$\DD_{u_1}\cdots\DD_{u_s}(a)$ lies in $I(\pi)$.
This proves that $I(\pi)_\DD\bigcap Z_\eps = m(\chi)_\DD$
and
$Z_\eps/m(\chi)_\DD= Z_\eps/(I(\pi)_\DD\bigcap Z_\eps)$.
Denote 
$$R_{\chi,\DD}=\frac{R_\eps}{m(\chi)_\DD R_\eps}.$$
{\bf Lemma 5.1}. Let $f$ be a nonzero element of $\FF$.
There exists $p\in\Omega_\chi$ 
such that  $f(p)\ne 0$ and $ R_\eps/m(\chi)R_\eps$ is isomorphic to
$ R_\eps/m(p)R_\eps$.\\
{\bf Proof}. 
If $f(\chi)\ne 0$, we put $p=\chi$. Let $f(\chi)=0$.
Consider $\Rpd$ as a $\FF$-module.
Let $W_1$ be the subvariety of points of $\MM_\chi$ for which the rank
of $ \Rcd$ over $\FF$ is non-minimal.
The set $W_1$ is a Poisson closed subset (i.e. an annihilator of
Poisson ideal) of $\MM_\chi$ [BD,4.1]. 
Since $m(\chi)_\DD$ is the minimal Poisson ideal which contains $m(\chi)$,
the point $\chi$ does not belong to $W_1$.
Since $m(p)_\DD=m(\chi)_\DD$ for all $p\in\Omega_\chi$ [BD,3.5],
then $\Omega_\chi\bigcap W_1=\emptyset$.
Denote $W_2=(\MM_\chi)_{\rm{sing}}$. 
As above $W_2$ is a closed Poisson
subvariety (Section 4) and  $\Omega_\chi\bigcap W_2=\emptyset$.
The $\FF$-module $\Rcd$ has constant rank over
the open smooth subset $\OO=\MM_\chi-(W_1\bigcup W_2)$.
It implies that $\Rcd$ is locally free over $\OO$
[H, Exercise II.5.8]. Let $U$  be an open subset such that
$\chi\in U\subset\OO$ and 
$$\Rcd'=\Rcd\otimes_{\FF}\Cb[U]$$ is free over $\Cb[U]$.
We consider
the algebra $C^{\infty}[U]$ of complex valued smooth functions 
on $U$ and extend $\Rcd'$ to the algebra
$$\Rcd^\infty=\Rcd'\otimes_{\Cb[U]}\Cb^\infty[U].$$
We consider this algebra as an algebra of smooth sections of 
a vector bundle with fibres $R_\eps/m(x)R_\eps$, $x\in U$.
For any $u\in R_\eps$ with $u_\eps=u\bmod(q-\eps)\in Z_\eps$ we extend
$\DD_u$ to the derivation of $\Rcd^\infty$. There exists a local flow 
$G_u(t)$ of automorphisms of $\Rcd^\infty$ 
lifting $\DD_u$  [DC-L,9.1],[DC-P,11.8],[BG,4.2].
The local flow is defined for small
$\vert t\vert$ and sends fibres to fibres.
Since $f$ is a nonzero element of $\FF=Z_\eps/{m(\chi)_\DD}$, 
there exist elements $u_1,\ldots,u_s$ such that 
 $u_{1\eps},\ldots, u_{n\eps}\in Z_\eps$ and  
the element
$\DD_{u_1}\cdots\DD_{u_s}(f)(\chi)\ne 0$. 
Let $G_{u_1}(t_1),\ldots, G_{u_s}(t_s)$
be the corresponding local flows defined for small $\vert t_i\vert$. 
Put
$p=G_{u_s}(t_s)\cdots G_{u_1}(t_1)(\chi)$.
We have $f(p)\ne 0$.
The algebras $R_\eps/m(\chi)R_\eps$ and $R_\eps/m(p)R_\eps$ are isomorphic. 
$\Box$\\ 
{\bf Definition 5.2}. 
We say that two irredicible representations $\pi_1$ and $\pi_2$ are
$\DD$-equivalent if $I(\pi_1)_\DD= I(\pi_2)_\DD$.\\
{\bf Theorem 5.3}. As above $R$ is a 
quantum solvable algebra obeying Conditions
3.2-3.4 and $\eps$ is an admissible $l^{\rm{th}}$ root of unity. We assert that:\\
1) dimension of any irreducible representation $\pi$ of $R_\eps$ equals
to $l^{\frac{1}{2}\mbox{dim}(\Omega_\chi)}$ where $\chi:Z_\eps\to\Cb$ is 
the central character of $\pi$ and $\Omega_\chi$ is the symplectic leaf of
$\chi$ in $M=\rm{Maxspec}Z_\eps$;\\  
2) any symplectic
leaf $\Omega_\chi$ is Zariski-open in its Zariski closure
(in the paper [BD] such leaves are called algebraic);
\\
3) if two points $\chi_1$ and $\chi_2$ belong to a common symplectic leaf,
then the algebras $R_{\chi_1}$ and $R_{\chi_2}$
are isomorphic;\\
4) let $\pi,\pi'$ are irreducible representations of $R_\eps$ and
$\chi,\chi'$ their central characters.
If $\pi$ and $\pi$ are $\DD$-equivalent, 
then $\chi$ and $\chi'$ lie in the 
common symplectic leaf.\\
{\bf Proof}. Our first goal is to prove 1).
Let $\pi$ be an irreducible representation of $R_\eps$.  
The ideal $I(\pi)_\DD$ is completely prime (Corollary 3.25).
There exists a unique standart 
semiprime $\DD$-stable ideal $\Pme$ and a denominator
set $\Sme$ in $\Rme=R_{\eps}/\Pme$ such that\\
1) $\Bme=\Rme\Sme^{-1}$ is a factor algebra of an
algebra twisted Laurent polynomials;\\
2) $\Ipd\supseteq\PP_{\mu\eps}$ and 
$\Sme\bigcap I/\Pme=\emptyset$.
(see Theorem 3.6, Definitions 3.10, 3.20).
We assume that $\Sme$ lies in $\Zme=\rm{Center}(\Rme)$ (see Remark 3.7).
We have
$$
\Rme\supset\Zme\supset \Zme'=\frac{Z_\eps}{\Pme\bigcap Z_\eps}.
\eqno(5.1)
$$
The algebra $\Zme$ has a Poisson structure and the algebra 
$\Zme'$ is a Poisson subalgebra of 
$\Zme$ (Lemma 3.12). 
The algebra $\Zme$ provides the quantum adjoint
action on $\Rme$. We shall denote it by $\oDD_u$ where $u\in R/\PP_\mu$ and
$u_\eps\in\Zme$.

By Lemma 3.16, 
$\Bme=A_\eps\otimes Z(B_\mu)_\eps$ where $A_\eps$ is an algebra of twisted 
Laurent polynomials generated by
$y_{1\eps}^{\pm 1},\ldots,y_{2r\eps}^{\pm 1}$ with 
relations (3.8). The $Z(B_\mu)_\eps$ is the reduction of $Z(B_\mu)$
modulo $q-\eps$. As in proof of Lemma 3.16, $Z(B_\mu)_\eps$  is generated by 
$z_{1\eps}^\pm 1,\ldots,z_{t\eps}$.
 The center $Z(\Bme)$ of $\Bme$ coincides with $\Zme\Sme^{-1}$
 and is generated by $y_{1\eps}^{\pm l},\ldots,y_{2r\eps}^{\pm l}$ and 
 $Z(B_\mu)_\eps$. We preserve the notation $\oDD$ for the
quantum adjoint action of $Z(\Bme)$ in $\Bme$ (Section 3) .\\
{\bf Step 1}.
By the Schur Lemma, $\pi(s)\in\Cb$ for all $s\in\Sme$.
Our goal in Step 1 is to reduce the general case to the case
$\pi(\Sme)\ne 0$. 

The ideal $J=\Ipd/\Pme$ of $\Rme$ is the greatest $\DD$-stable ideal 
in $I(\pi)/\Pme$. The ideal $J$ is is completely prime.
The localization $J\Sme^{-1}$ is $\DD$-stable ideal in $\Bme$.  
By Corollary 3.22, $J\Sme^{-1}$ is $\oDD$-stable.
This implies that $J$ is also $\oDD$-stable with respect to
the quantum adjoint action of $\Zme$ in $\Rme$.
Denote
$$
\Rpd=R_\eps/\Ipd=\Rme/J , \qquad \ZZ=\frac{\Zme}{J\bigcap \Zme}.
$$
We preserve the notation $\Sme$ for the its image in $\ZZ$ and denote
$\Bpd=\Rpd\Sme^{-1}$.
Remark that
$$
\ZZ\Sme^{-1}=\frac{\Zme\Sme^{-1}}{J\Sme^{-1}\bigcap \Zme\Sme^{-1}}=
\frac{Z(\Bme)}{J\Sme^{-1}\bigcap Z(\Bme)},\eqno (5.2)$$
$$
\frac{\Zme'}{J\bigcap \Zme'}=\frac{Z_\eps}{\Ipd\bigcap Z_\eps}=
\frac{Z_\eps}{m(\chi)_\DD}=\Cb[M_\chi]=\FF.
$$
We have
$$
\Rpd\supset\ZZ\supset\FF.\eqno(5.3)
$$
All algebras in (5.3) are domains, $\ZZ$ is a Poisson algebra, 
$\FF$ is its Poisson subalgebra.
Let $s_1,\ldots, s_k$ be $q$-commuting generators
of $\Sme$.
Since $\ZZ$ is finite over $\FF$, the $\ZZ \FF^{-1}$ is a field.
The element $s_{i}^{-1}$ can be presented in the form
$$s_{i}^{-1}=a_iw_i^{-1},
$$
where $a_i\in \ZZ$, $w_i\in \FF$.
Let $N$ denote the denominator subset in $\FF$ 
generated by $w_1,\ldots,w_k$.
After localization we obtain
$$\Bpd N^{-1}=\Rpd N^{-1}\supset 
\Aa:=\ZZ N^{-1}\supset \Bb:=\FF N^{-1}.\eqno(5.4)
$$
Denote $X=\Ms(\Aa)$, $Y=\Ms(\Bb)$ and
$f_0=w_1\cdots w_k\in \FF.$
We apply Lemma 5.1. 
To prove statement 1) of the Theorem we may consider that $f_0(\chi)\ne 0$.
Then $\pi(w_i)=w_i(\chi)\ne 0$. Since $w_i=a_is_i$, $\pi(s_i)\ne 0$.\\
{\bf Step 2}. By Step 1, one may suppose that $\pi$ is lifted to a
representation of $\Bpd N^{-1}$. 
Therefore, $\pi$ is an irreducible representation of the algebra
of twisted Laurent polynomials. By Lemma 3.16,
$$\dim \pi= l^r.\eqno (5.5)$$
We denote $\pi(z_{i\eps})=\alpha_i\in\Cb$.
According to Proposition 3.21, an ideal in
$\Bme N^{-1}$ is  $\oDD$-stable ideal if and only if it is
 $\DD$-stable.
 This implies that $JN^{-1}$ is the greatest $\oDD$ stable ideal in
 $(I(\pi)/\Pme)N^{-1}$
 and
$$
 JN^{-1}=\sum_{i=1}^t \Bme N^{-1}(z_i-\alpha_i).
 $$
This implies that the Poisson algebra $\Aa$ is the localization
of subalgebra generated by 
$y_{l\eps}^{\pm 1},\ldots,y_{2r\eps}^{\pm l}$.
The set $X$ is a Zariski-open subset
of $\Cb^{2r}$. By Lemma.3.8, $X$ is a symplectic variety.
 We have $2r=\dim X= \rk \Aa=\rk_x\Aa$ for any $x\in X$.
 
 The set $Y$ is a Zariski-open subset in $M_\chi=\overline{\Omega_\chi}$.  
 Since $\Aa$ is finite over $\Bb$, then $\dim X=\dim Y$ and,
  by Proposition 4.2, 
 there exists a Zariski-open subset $U$ of $Y$ (we may  assume that 
 $U\bigcap (M_\chi)_{\rm{sing}}=\emptyset$) such that
 $\rk_y\Bb=\rk_x\Aa=\dim X$ for $y=\phi(x)$.
 We obtain $\rk_y\Bb=\dim Y=\dim U$ for all $y\in U$.
 Hence, $ U$ is contained in a symplectic leaf.
 Since $\Omega_\chi\bigcap U\ne \emptyset$, then $\Omega_\chi\supset U$
 and $\dim\Omega_\chi=\dim U=\dim Y=\dim X$.
 Finally, we get 
 $r=\frac{1}{2}\dim X=\frac{1}{2}\dim \Omega_\chi$ and, by (5.5),
$\dim \pi=l^{\frac{1}{2}\dim{\Omega_\chi}}$.This proves 1).\\
{\bf Step 3}. Our goal is to prove 2),3),4).\\
2) We have proved in Step 2 that $\Omega_\chi$ contains a Zariski-open
subset $U$. For any element  $p\in\overline{\Omega_\chi}-\Omega_\chi$
we have
$$
\dim\Omega_p\le \dim(\overline{\Omega_\chi}-U) < \dim U=\dim \Omega_\chi.$$
The set $\MM_{<2r}:=\{p\in\MM: \mbox{dim}(\Omega_p)< 2r\}$ 
is Zariski-closed in $\MM$ [BG,3.1].
The subset $\overline{\Omega_\chi}-\Omega_\chi$  
coincides with
$\overline{\Omega_\chi}\bigcap \MM{<2r}$ and is Zariski-closed.
This proves statement 2.\\
3) For any two points $\chi_1$ and $\chi_2$ of $\Omega_\chi$. 
By Lemma 5.1 , there exist $p_1,p_2\in U$ such that
$R_{\chi_1}$ is isomorphic to $R_{p_1}$
and $R_{\chi_2}$ is isomorphic to $R_{p_2}$.
The algebras $R_{p_1}$ and $R_{p_2}$ are isomorphic as fibers of
the algebra
of twisted Laurent polynomials $\Bpd N^{-1}$.  This proves 3.\\
4) Let $\pi_1$ and $\pi_2$ be two $\DD$-equivalent representations of
$R_\eps$. Then $I(\pi_1)_\DD=I(\pi_2)_\DD$. Then
$m(\chi_1)_\DD=m(\chi_2)_\DD$. By 2), the points $\chi_1$ and $\chi_2$ lie in 
a common symplectic leaf.$\Box$

\end{document}